\documentclass{article}

\newif\ifarxiv
\arxivtrue

\usepackage{graphicx}
\usepackage{amsmath,amssymb,amsthm,bm}
\usepackage[a4paper]{geometry}
\usepackage{algorithm,algorithmic}
\usepackage[colorlinks=false,pdfborder={0 0 0}]{hyperref}
\usepackage{authblk}
\usepackage{color}
\usepackage{subcaption}

\graphicspath{{fig/}}
\DeclareGraphicsExtensions{.mps,.pdf,.eps,.png}

\newenvironment{keywords}{\medskip\textbf{Keywords:}}{}
\newenvironment{AMS}{\medskip\textbf{AMS subject classifications (2020).}}{}

\newtheorem{theorem}{Theorem}

\newtheorem{remark}{Remark}
\newtheorem{lemma}{Lemma}

\newtheorem{corollary}{Corollary}
\theoremstyle{plain}

\DeclareMathOperator*{\argmin}{\arg\min}

\DeclareMathOperator{\diag}{diag}

\DeclareMathOperator{\Span}{span}

\newcommand*{\set}[1]{\left\lbrace#1\right\rbrace}

\newcommand*{\trans}{^{\top}}

\newcommand*{\itrans}{^{-\top}}

\newcommand{\bmat}[1]{\begin{bmatrix}#1\end{bmatrix}}

\newcommand*{\abs}[1]{\bigl\lvert#1\bigr\rvert}
\newcommand*{\absbig}[1]{\biggl\lvert#1\biggr\rvert}
\newcommand*{\norm}[1]{\bigl\Vert#1\bigr\rVert}
\newcommand*{\normbig}[1]{\biggl\Vert#1\biggr\rVert}
\newcommand*{\normA}[1]{\bigl\Vert#1\bigr\rVert_A}
\newcommand*{\normMinv}[1]{\bigl\Vert#1\bigr\rVert_{M^{-1}}}

\def\adots{\mathinner{\mkern2mu\raise1pt\hbox{.}\mkern2mu
    \raise4pt\hbox{.}\mkern2mu\raise7pt\hbox{.}\mkern1mu}}

\newcommand{\lamdmax}{\lambda_{\max}}
\newcommand{\lamdmin}{\lambda_{\min}}

\newcommand*{\macheps}{\bm u}
\newcommand*{\machepss}{\bm u_s}
\newcommand*{\machepsq}{\bm u_q}

\newcommand*{\machepsleft}{\bm u_L}
\newcommand*{\machepsright}{\bm u_R}

\newcommand*{\epss}{\varepsilon_{\mathrm{s}}}
\newcommand*{\epsq}{\varepsilon_{\mathrm{q}}}
\newcommand*{\epsz}{\varepsilon_{\mathrm{z}}}
\newcommand*{\epspre}{\epsilon_{\mathrm{pre}}^{(s,q)}}
\newcommand*{\epspreq}{\epsilon_{\mathrm{pre}}^{(s,q)}}

\newcommand*{\bigO}{O}
\DeclareMathOperator{\fl}{f{}l}

\newcommand*{\Dalpha}{\delta{\alpha}}
\newcommand*{\Dhatalpha}{\delta{\hat{\alpha}}}
\newcommand*{\Dx}{\Delta{x}}
\newcommand*{\Dbarx}{\Delta{\bar{x}}}
\newcommand*{\Dr}{\Delta{r}}

\newcommand*{\Ds}{\Delta{s}}
\newcommand*{\Dq}{\Delta{q}}
\newcommand*{\Dbeta}{\delta{\beta}}
\newcommand*{\Dhatbeta}{\delta{\hat{\beta}}}
\newcommand*{\Dp}{\Delta{p}}
\newcommand*{\Drp}{\delta(r\trans p)}

\newcommand*{\krylov}{\mathcal{K}}
\newcommand*{\kstar}{k^{\star}}

\newcommand*{\hatalpha}{\hat{\alpha}}
\newcommand*{\barx}{\bar{x}}
\newcommand*{\barr}{\bar{r}}
\newcommand*{\hatx}{\hat{x}}
\newcommand*{\hatr}{\hat{r}}

\newcommand*{\hatq}{\hat{q}}
\newcommand*{\hatbeta}{\hat{\beta}}
\newcommand*{\hatp}{\hat{p}}

\newcommand*{\barxkone}{\bar{x}_{k+1}}
\newcommand*{\barrkone}{\bar{r}_{k+1}}
\newcommand*{\hatxkone}{\hat{x}_{k+1}}
\newcommand*{\hatrkone}{\hat{r}_{k+1}}

\newcommand*{\hatskone}{\hat{s}_{k+1}}
\newcommand*{\hatqkone}{\hat{q}_{k+1}}
\newcommand*{\hatbetakone}{\hat{\beta}_{k+1}}
\newcommand*{\hatpkone}{\hat{p}_{k+1}}
\newcommand*{\hatalphak}{\hat{\alpha}_{k}}
\newcommand*{\barxk}{\bar{x}_{k}}
\newcommand*{\barrk}{\bar{r}_{k}}
\newcommand*{\hatxk}{\hat{x}_{k}}
\newcommand*{\hatrk}{\hat{r}_{k}}

\newcommand*{\hatqk}{\hat{q}_{k}}

\newcommand*{\hatpk}{\hat{p}_{k}}

\title{Forward and backward error bounds for a mixed precision preconditioned conjugate gradient algorithm}
\author[1]{Thomas Bake}
\author[1]{Erin Carson}
\author[1]{Yuxin Ma}
\affil[1]{Department of Numerical Mathematics, Faculty of Mathematics and Physics, Charles University, Sokolovsk\'{a} 49/83, 186 75 Praha 8, Czechia

{\tt Email: bake@karlin.mff.cuni.cz, carson@karlin.mff.cuni.cz, yuxin.ma@matfyz.cuni.cz}}

\begin{document}
\maketitle

\begin{abstract}
The preconditioned conjugate gradient (PCG) algorithm is one of the most popular algorithms for solving large-scale linear systems \(Ax = b\), where \(A\) is a symmetric positive definite matrix.
Rather than computing residuals directly, it updates the residual vectors recursively.
Current analyses of the conjugate gradient (CG) algorithm in finite precision typically assume that the norm of the recursively updated residual goes orders of magnitude below the machine precision, focusing mainly on bounding the residual gap thereafter.
This work introduces a framework for the PCG algorithm and provides rigorous proofs that the relative backward and forward errors of the computed results of PCG can reach the levels \(\bigO(\macheps)\) and \(\bigO(\macheps)\kappa(A)^{1/2}\), respectively, after a sufficient number of iterations without relying on an assumption concerning the norm of the recursively updated residual, where \(\macheps\) represents the unit roundoff and \(\kappa(A)\) is the condition number of \(A\).
Our PCG framework further shows that applying preconditioners in low precision does not compromise the accuracy of the final results, provided that reasonable conditions are satisfied.
Moreover, this framework introduces a new split PCG variant that improves upon the classical split PCG algorithm when the left preconditioner is applied in low precision.
Our theoretical results are illustrated through a set of numerical experiments.

\begin{keywords}
conjugated gradient algorithm, forward and backward error analysis, finite precision arithmetic, mixed precision
\end{keywords}

\begin{AMS}
65F10, 65F08, 65G50, 65Y20
\end{AMS}
\end{abstract}

\section{Introduction}
\label{sec:introduction}
Given a symmetric positive definite (SPD) matrix \(A \in \mathbb R^{n\times n}\) and a right-hand side \(b\in\mathbb R^n\), we study the approximation to the solution of linear systems
\begin{equation} \label{problem:linear-sys}
    Ax=b.
\end{equation}
In particular, when the matrix $A$ is large and sparse, the conjugate gradient (CG) algorithm is the method of choice. 
The most widely used CG variant was introduced by Hestenes and Stiefel~\cite{HS1952}.
Although there are numerous mathematically equivalent versions of the CG algorithm, such as the ones described in~\cite{GS2000}, the implementation of Hestenes and Stiefel~\cite{HS1952}, detailed in Algorithm~\ref{alg:cg}, is still the most popular because of its simplicity and its superior behavior in finite precision.

\begin{algorithm}[H]
\begin{algorithmic}[1]
    \caption{Conjugate Gradient (CG) algorithm \label{alg:cg}}
    \REQUIRE
     A symmetric positive definite matrix \(A \in \mathbb R^{n\times n}\), a right-hand side \(b \in \mathbb R^{n}\), an initial guess $x_0 \in \mathbb{R}^{n}$, and maximum number of iterations \(\texttt{maxiter}\).
    \ENSURE
    Approximate solution \(x_{k + 1} \in \mathbb R^{n}\) of \(A x = b\).

    \STATE \(r_0 = b - Ax_0\).
    \STATE \(p_0 = r_0\).
    \FOR{\(k = 0, 1, 2, \dotsc, \texttt{maxiter}\)}
        \STATE \(\alpha_{k} = \frac{r_{k}\trans r_{k}}{p_{k}\trans Ap_{k}}\).
        \STATE \(x_{k+1} = x_{k} + \alpha_{k}p_{k}\).
        \STATE \(r_{k+1} = r_{k}-\alpha_{k}Ap_{k}\). \label{line:updater}
        \IF{the stopping criterion is satisfied}
            \RETURN \( x_{k+1}\).
        \ENDIF
        \STATE \(\beta_{k+1} = \frac{r_{k+1}\trans r_{k+1}}{r_{k}\trans r_{k}}\).
        \STATE \(p_{k+1} = r_{k+1} + \beta_{k+1}p_{k}\).
    \ENDFOR
\end{algorithmic}
\end{algorithm}

It is well known that the application of a preconditioner \(M \in \mathbb R^{n \times n}\) has the potential to accelerate the convergence of iterative methods, and their application is indispensable in practical computations. In the context of CG, the effect of preconditioning is to transform \eqref{problem:linear-sys} to the equivalent linear system
\begin{equation}\label{problem:prec-linear-sys}
M_L^{-1}AM_R^{-1} y = M_L^{-1}b, \quad \text{with} \quad M_R^{-1} y = x, 
\end{equation}
where the matrices \(M_L\) and \(M_R \in \mathbb R^{n \times n}\), $M_L M_R=M$, are called left and right preconditioners, respectively.
In practice, the preconditioned system \eqref{problem:prec-linear-sys} is not explicitly formed. Instead, the iterative method is reformulated so that at every iteration we apply the preconditioners to one or more vectors.
Due to the practical effectiveness of the preconditioning technique, we will focus on the preconditioned conjugate gradient (PCG) algorithm in this work.

Rounding errors, as is widely known, can lead to an insufficiently accurate solution and slow down the convergence process for CG~\cite{GS2000}.
Consequently, conducting a rounding error analysis is crucial to understand the behavior of CG.
In general, CG does not directly update \(r_k\) using \(b-Ax_k\).
Specifically, in Line~\ref{line:updater} of Algorithm~\ref{alg:cg}, \(r_k\) is revised through a method that is mathematically equivalent but numerically distinct, resulting in a gap between the (computed) \emph{true residual} \(\tilde{r}_k = b - A \hat{x}_k\) and the (computed) \emph{recursively updated residual} \(\hat{r}_k\), where $\hat{\cdot}$ denotes the computed quantity.
Existing accuracy analyses of CG in finite precision assume that the norm of the recursively updated
residual can eventually reach the level \(\bigO(\macheps)\) after subsequent iterations, where \(\macheps\) denotes the unit roundoff.
This assumption is commonly used to determine the maximal attainable accuracy in CG-like algorithms together with residual gap bounds; see~\cite{Greenbaum1997,GS2000}.

In this work, we present a framework for the PCG algorithm based on Algorithm~\ref{alg:cg}, and investigate its behavior in finite precision arithmetic.
It is important to note that our analysis is not based on an assumption on \(\norm{\hatrk}\), where \(\norm{\cdot}\) denotes the \(2\)-norm of a vector or the spectral norm of a matrix.
Instead, we use the well-known fact that CG can be seen as an optimization method that minimizes a quadratic function,
\begin{equation}\label{eq:quadratic-function}
     f(y) := \frac{1}{2} \langle y, y \rangle_A - \langle y, b \rangle = \frac{1}{2} y^\top A y - y^\top b.
\end{equation}
More precisely, we show that there exists a \(\kstar\) such that the difference of the corresponding function values in iterations \(\kstar\) and \(\kstar+1\) can be sufficiently small.
In this specific iteration \(\kstar\), we first show that, under reasonable assumptions, \(\norm{\hatr_{\kstar}}/\bigl(\norm{A}\norm{x}\bigr)\) can be bounded by \(\bigO\bigl(n(\kstar)^2\macheps\bigr)\kappa(M)^{1/2}\), where \(\kappa(M):=\norm{M}\cdot \norm{M^{-1}}\) is the condition number of the matrix \(M\) with respect to the \(2\)-norm.
Furthermore, we show that, under reasonable assumptions, unpreconditioned CG, i.e., Algorithm~\ref{alg:cg}, can produce a relative backward error, i.e., \(\norm{\tilde{r}_{\kstar}}/\bigl(\norm{A}\norm{x}\bigr)\), that is bounded by \(\bigO\bigl(n(\kstar)^2\macheps\bigr)\).
If, in addition, we consider using a preconditioner $M$ in our PCG framework, then the bound has the additional factor of $\kappa(M)^{1/2}$.
Moreover, we prove that for the same iteration, a comparable level can be reached by the $A$-norm of the relative forward error, which can also be formulated in the \(2\)-norm.
The results mentioned above allow us to derive conditions on the precisions with which the preconditioners have to be applied such that forward and backward errors of the same magnitude are produced in a certain iteration.
Thus, we can theoretically guarantee that we can employ low precision for the application of the preconditioners while retaining working precision accuracy for the computed approximation, as long as reasonable conditions are satisfied.

Due to the development of hardware and machine learning techniques, mixed precision algorithms have attracted significant attention in recent years. 
Rather than the traditional approach of using a uniform precision (the \emph{working precision}) for storage and computation throughout an algorithm, it is desirable to use low precision in select parts of a computation, as this can require less memory, less data movement, less energy, and less computational time.
Since using low precision may compromise the accuracy of the final result in finite precision, carefully adopting low precision arithmetic is essential to significantly accelerate computation while maintaining the required accuracy; see the surveys~\cite{Survey2021,HM2022} for more information.

In PCG, the application of preconditioners can involve, for example, solving two triangular systems when the preconditioner is given by an incomplete Cholesky factorization. Therefore, if $n$ is very large and the preconditioners are dense relative to $A$, the cost of applying the preconditioners in each iteration can be significant. Following this example, if the generated preconditioner contains significant fill-in---the introduction of nonzero entries in the Cholesky factors, whose corresponding entries in the original matrix are zero---then its storage might not even be feasible for large-scale problems. These drawbacks practically motivate the storage and application of preconditioners in low precision.
Our numerical experiments primarily focus on using incomplete Cholesky factorization for preconditioner computation, a process that can be reliably performed in low precision~\cite{ST20224}.

The remainder of this paper is organized as follows.
In Section~\ref{sec:background}, we introduce some background and notation, and give an overview of related work.  
Later, in Section~\ref{sec:algorithm}, we present a comprehensive framework for the preconditioned CG algorithm based on Algorithm~\ref{alg:cg}.
This includes implementations with left, right, and split preconditioners. Then Section~\ref{sec:stability} offers a backward and forward error analysis of the preconditioned CG algorithm.
Numerical experiments are presented in Section~\ref{sec:experiments} to confirm our theoretical results.
\section{Notation, background and related work}\label{sec:background}

Throughout this paper, we use the fact that symmetric positive definite matrices induce an inner product in $\mathbb R^n$ with corresponding norm, and that the 2-norm condition number is given by the ratio of its largest to smallest eigenvalue.
That is, for an SPD matrix $A \in \mathbb R^{n \times n}$ we have the \(A\)-inner product and the \(A\)-norm:
\[
    \langle x, y\rangle_A := y^\top A x, \qquad \norm{x}_A = \langle x, x \rangle_A^{\frac{1}{2}} = \norm{A^{\frac{1}{2}}x}, \qquad \text{for all } x, y \in \mathbb R^n, 
\]
and \( \kappa(A) = \norm{A} \norm{A^{-1}} = \lambda_{\max} (A) / \lambda_{\min} (A) \).
In addition, we use $\hat{\cdot}$ to denote computed quantities. The unit roundoff corresponding to working precision is denoted by $\macheps$. 
Moreover, for the unit roundoff corresponding to the application of the preconditioners, we add a subscript.
For example, the vector $\hat{s}_k = M^{-1}_L \hatrk + \Delta s_k$ is computed with the precision associated with the unit roundoff $\macheps_s$.

Given an initial guess $x_0$, at each iteration $k$, the CG algorithm in exact arithmetic computes an approximate solution $x_k \in x_0 + \krylov_k (A, r_0)$ of $x$ in \eqref{problem:linear-sys}, and is mathematically characterized by the equivalent optimality and orthogonality properties
\begin{equation}\label{eq:cgcharac}
  \norm{x - x_k}_A = \min_{z \in x_0 + \krylov_k (A, r_0)} \norm{x - z}_A 
  \quad \text{and} \qquad
  b - A x_k \perp \krylov_k (A, r_0),
\end{equation}
where
\begin{equation*}
  \krylov_k (A, r_0) := \Span\{r_0, A r_0, A^2 r_0, \dots, A^{k - 1} r_0\}
\end{equation*}
is the $k$-dimensional Krylov subspace of $A$ and $r_0$.
From the optimality characterization of the $k$-th CG iterate in \eqref{eq:cgcharac}, one can easily relate the error in the $A$-norm to the polynomial minimization problem
\begin{equation}\label{eq:polynomial-min}
  \norm{x - x_k}_A = \min_{p \in \mathcal{P}_k(0)} \norm{p(A)(x- x_0)}_A,
\end{equation}
where $\mathcal{P}_k(0)$ is the space of $k$-th degree polynomials that are one at the origin. If $A = Y \Lambda Y \trans$ is an orthogonal diagonalization of the SPD matrix $A$ with $\Lambda = \diag(\lambda_1, \dots, \lambda_n)$, and $\eta_1, \dots, \eta_n$ are the components of $r_0$ in the corresponding invariant subspaces (that is, $r_0 = Y [\eta_1, \dots, \eta_n]\trans$), then one can express the right-hand side of~\eqref{eq:polynomial-min}---and hence the error in the $A$-norm---as 
\begin{equation}\label{eq:polynomial-eigvals-min}
	\normA{x - x_k} = \min_{p \in \mathcal{P}_k(0)} \left( \sum_{i = 1}^n \eta_i^2 \frac{p(\lambda_i)^2}{\lambda_i} \right)^{1/2};
\end{equation}
see~\cite[Section 2]{CLS2024}.
The polynomial that minimizes the above sum is called the \emph{CG polynomial} and it is denoted by $\varphi^{\text{CG}}_k$, and its roots are called the \emph{Ritz values} of $A$\footnote{These are also the eigenvalues of the tridiagonal matrix that is (implicitly) generated by the CG algorithm and the Lanczos process at each iteration.}.
This expression will be helpful to explain the experimental setting in Section~\ref{sec:experiments}.

It is well known that the CG algorithm has multiple mathematical interpretations.
It can be seen as a projection process onto Krylov subspaces, as the Lanczos method for computing eigenvalues, as the Gauss--Christoffel quadrature of a Riemann--Stieltjes integral with respect to a certain distribution function with given points of increase, as the Stieltjes recurrence for generating orthonormal polynomials, and more.
For a detailed exploration of these interconnections, see~\cite[Section 3]{LS2012}. 

In exact arithmetic, the recursively updated residuals \(r_0\), \(r_1\), \(\dotsc\), \(r_k\) and the direction vectors \(p_0\), \(p_1\), \(\dotsc\), \(p_k\) form an orthonormal basis with respect to the Euclidean and $A$-inner products of $\krylov_k (A, r_0)$, respectively.
In addition, Krylov subspace methods, such as CG, can be considered mathematically as direct methods, since they are guaranteed to terminate in at most $n$ iterations.
However, in finite precision, rounding errors cause loss of orthogonality among the vectors that (ideally) span the Krylov subspace as the iterations proceed, so that their attractive properties in exact arithmetic do not necessarily hold anymore. 
Therefore, a great deal of effort has been put into investigating the numerical stability of Krylov subspace methods through their \emph{forward} and \emph{backward errors}; see, e.g., \cite{JR2008,PRS2006,PS2002}.
In particular, the HSCG algorithm is not guaranteed to produce a backward stable solution in $n$ iterations, in contrast to full recurrence Krylov subspace methods like MGS-GMRES, which was proved to be backward stable in~\cite{PRS2006}.
Generally speaking, the use of short recurrences in HSCG may cause orthogonality to be completely lost, resulting in the rank deficiency of the computed basis of the Krylov subspaces.
This in turn can cause a dramatic delay of convergence in the CG algorithm.
In~\cite{Greenbaum1989}, Greenbaum studied this phenomenon and gave a backward-like analysis of CG (and/or Lanczos). 
In her backward-like analysis, she showed that a perturbed CG algorithm applied to a given problem with an SPD matrix $A \in \mathbb R^{n \times n}$ and an initial guess $x_0 \in \mathbb R^n$ is equivalent---up to a negligible difference---to the exact CG algorithm applied to a larger problem.
More precisely, in the larger problem, the matrix has its eigenvalues tightly clustered around the eigenvalues of $A$; see \cite[Section 2.3]{CLS2024} for an example of such a construction.
Furthermore, she used the fact that both the Lanczos and CG algorithms (implicitly) generate orthogonal polynomials whose roots (Ritz values) progressively approximate the eigenvalues of the given SPD matrix.
In the interest of space, we will not further elaborate on Greenbaum's creative and thorough work, so we point the reader to her original paper and to \cite[Section 4.3]{MS2006} for a comprehensive summary of her work and an excellent survey of the CG and Lanczos algorithms in exact and finite precision arithmetic.

Another complication caused by the accumulation of rounding errors is the potential decrease of the attainable accuracy of CG-like algorithms.
In~\cite{Greenbaum1997}, Greenbaum showed that for iterative methods that update their iterates and residuals using recurrences of the form
\[
    x_{k + 1} = x_k + \alpha_k p_k \quad \text{ and } \quad r_{k + 1} = r_k - \alpha_k A p_k,
\]
the decrease of the true residuals in norm is limited by the largest ratio of the norm of an iterate to that of the true solution $x$ of \eqref{problem:linear-sys}, i.e., \(\max_{j\leq k}\bigl(\norm{\hatx_j}/\norm{x}\bigr)\).
She assumed that, after some number of iterations, the norm of the recursively updated residual $\hatrk$ goes levels below the machine precision and estimated the norm of the residual gap $\tilde{r}_k - \hat{r}_k$\footnote{In finite precision, it has been well established that CG-like algorithms might reach a point at which $\norm{b - A \hatxk}$ and $\norm{\hatrk}$ start to deviate from each other. More precisely, it has been observed in many practical settings that $\norm{b - A \hatxk}$ stagnates while $\norm{\hatrk}$ continues to decrease several levels below the machine precision.}, where the true residual \(\tilde{r}_k\) and the recursively updated residual $\hatrk$ are defined in Section~\ref{sec:introduction}. Under this assumption, the norm of the residual gap determines the attainable accuracy since
\[
    \norm{b - A \hatxk} \leq \norm{b - A \hatxk - \hatrk} + \norm{\hatrk}.
\]

Numerous authors have used this approach (or analogous modifications depending on the problem at hand) to investigate the attainable accuracy of CG variants and other Krylov subspace methods; see, e.g., \cite{BES1998,GS2000,JR2008,MSV2001}.
More recently, the attainable accuracy has also been explored through this approach for CG variants developed for high-performance computing that reduce synchronization or communication; see, e.g., \cite{CD2014,CRSTT2018}.
In comparison to the delay of convergence of CG or other Krylov subspace methods, the attainable accuracy has not been given the same attention under the justified argument that in practical applications the process is stopped much earlier than the maximal attainable accuracy is reached (which is an advantage of  employing iterative over direct methods). 
However, when using low precision or parallel implementations of CG, the maximal attainable accuracy could be damaged to a point that it might reach a level that is orders of magnitude above the required tolerance of the practical application.
If this is the case, the problem might be numerically unsolvable for the given implementation of the method and the choice of precisions for the corresponding tasks.
Recently, the authors in~\cite{GSW2025} estimate the norm of the residual gap in an adaptive mixed precision PCG with a left preconditioner to alternate between double and single precision for the computation of matrix-vector products across iterations.

As mentioned above, our analysis uses one of the many representations of CG, namely that of an optimization method that minimizes the quadratic function~\eqref{eq:quadratic-function}.
Notice that the gradient of \(f\) is \(Ay-b\).
Therefore, the optimal solution \(x = \argmin_y f(y)\) also represents the solution to the linear system~\eqref{problem:linear-sys}.
From an optimization perspective, it can be assured that the function \(f\) will converge to its minimum.
Equivalently, we have \(f(x) = -\normA{x}^2/2\) and then
\begin{equation}\label{eq:energy-min}
\begin{split}
    f(x_k) - f(x) 
    &= \frac{1}{2} \normA{x_k}^2 - \langle x_k, b \rangle + \frac{1}{2} \normA{x}^2 \\
    &= \frac{1}{2} \normA{x_k}^2 - \langle x_k, x \rangle_A + \frac{1}{2} \normA{x}^2 \\
    &= \frac{1}{2} \normA{x - x_k}^2,
\end{split}
\end{equation}
so that minimizing the $A$-norm of the error is equivalent to minimizing the quadratic function.

\section{A framework for the preconditioned CG algorithm}
\label{sec:algorithm}
In~\cite[Section 9]{Saad2003}, Saad discussed the left-, right-, and split-preconditioned CG (PCG) algorithms individually.
Here, we present a framework that accommodates all these PCG algorithms using various preconditioners, as depicted in Algorithm~\ref{alg:pcg}.
This framework aligns with the implementation of Hestenes and Stiefel and is designed to address the preconditioned linear system \eqref{problem:prec-linear-sys} which is equivalent to Algorithm~\ref{alg:cg} when no preconditioner is applied, i.e., when \(M_L=M_R=I\).
The CG algorithm is well defined for SPD matrices only. Hence, it is necessary to maintain symmetry in the preconditioned algorithm, which translates into enforcing that the matrices \(M^{-1} = M_R^{-1} M_L^{-1}\) and \(M_L^{-1}AM_R^{-1}\) are symmetric positive definite as well. Algorithm \ref{alg:pcg} achieves this since the inner products are symmetric with respect to the $M$-inner product $\langle x, y \rangle_M = y^\top M x$, where $x, y \in \mathbb R^n$. 
Within this framework, it becomes evident that PCG algorithms employing left and right preconditioners are mathematically and numerically equivalent, as noted in~\cite[Section 9]{Saad2003}. However, in contrast to the split preconditioned PCG algorithm presented by~\cite[Algorithm 9.2]{Saad2003}, our framework recovers both left and right PCG when setting the preconditioners accordingly. For instance, if we have a Cholesky factorization $M = L L^\top$, then setting 
\[
    M_L = L L^\top, M_R = I \quad \text{ or } \quad M_L = I, M_R = L L^\top
\] 
in Algorithm~\ref{alg:pcg} yields left or right PCG, respectively. This is not the case for~\cite[Algorithm 9.2]{Saad2003}. 

\begin{algorithm}[htbp!]
\begin{algorithmic}[1]
    \caption{Preconditioned Conjugate Gradient (PCG) algorithm \label{alg:pcg}}
    \REQUIRE
     A symmetric positive definite matrix \(A \in \mathbb R^{n\times n}\), a right-hand side \(b \in \mathbb R^{n}\), an initial guess $x_0 \in \mathbb{R}^{n}$, the left and right preconditioners \(M_L^{-1}\) and \(M_R^{-1}\) satisfying that \(M^{-1} = M_R^{-1} M_L^{-1}\) and \(M_L^{-1}AM_R^{-1}\) are symmetric positive definite matrices, and maximum number of iterations \(\texttt{maxiter}\).
    \ENSURE
    Approximate solution \(x_{k + 1} \in \mathbb R^{n}\) of \(A x = b\).

    \STATE \(r_0 = b - Ax_0\).
    \STATE \(s_0 = M_L^{-1} r_0\).
    \STATE \(q_0 = M_R^{-1} s_0\) and \(p_0 = q_0\).
    \FOR{\(k = 0, 1, 2, \dotsc, \texttt{maxiter}\)}
        \STATE \(\alpha_{k} = \frac{r_{k}\trans q_{k}}{p_{k}\trans Ap_{k}}\). \label{line:alpha}
        \STATE \(x_{k+1} = x_{k} + \alpha_{k}p_{k}\). \label{line:x}
        \STATE \(r_{k+1} = r_{k}-\alpha_{k}Ap_{k}\). \label{line:r}
        \STATE \(s_{k+1} = M_L^{-1} r_{k+1}\). \label{line:s}
        \STATE \(q_{k+1} = M_R^{-1} s_{k+1}\). \label{line:q}
        \IF{the stopping criterion is satisfied}
            \RETURN \(x_{k+1}\).
        \ENDIF
        \STATE \(\beta_{k+1} = \frac{r_{k+1}\trans q_{k+1}}{r_{k}\trans q_{k}}\). \label{line:beta}
        \STATE \(p_{k+1} = q_{k+1} + \beta_{k+1}p_{k}\). \label{line:p}
    \ENDFOR
\end{algorithmic}
\end{algorithm}

Furthermore, the split-preconditioned variant of Saad (see \cite[Algorithm 9.2]{Saad2003}) suffers from an amplification of the local rounding errors that arise from the computation of $\hatrk$ when the left preconditioner is applied in low precision. This significantly limits the accuracy of the corresponding approximation compared to the left and right preconditioned counterparts. The cause of this amplification is the \emph{recursive} application of the left preconditioner in the recurrence to generate the recursively updated residual $\hatrk$. In our framework, the application of the left preconditioner occurs in a subsequent step and therefore Algorithm~\ref{alg:pcg} avoids this problem.
This will be illustrated by numerical experiments provided in Section~\ref{sec:experiments}.

\section{The numerical stability of PCG}
\label{sec:stability}
In this section, we present a rounding error analysis of the PCG algorithm, i.e., Algorithm~\ref{alg:pcg}.
First, we perform a rounding error analysis for each step within Algorithm~\ref{alg:pcg} in Section~\ref{subsec:analysis-step}.
Subsequently, in Section~\ref{subsec:localorth}, we consider the local orthogonality between the computed results of \(r_{k+1}\) and \(p_k\).
In the final part, Section~\ref{subsec:back-forward-error} delves into the analysis of backward and forward errors associated with Algorithm~\ref{alg:pcg}, culminating in the summary of findings presented in Theorem~\ref{thm:backward-forward-err}.

In order to perform a rounding error analysis of PCG, we employ the standard model:
\begin{equation}
\label{eq:model0}
\fl(\alpha\circ\beta)=(\alpha\circ\beta)(1+\delta),
\qquad \lvert\delta\rvert\leq\macheps,
\end{equation}
for any finite floating-point numbers \(\alpha\), \(\beta\in\mathbb R\) and
\(\circ\in\set{+,-,\times,/}\).

Let us make some final notational remarks before proceeding.
In the following analysis, we assume \(n\macheps<1\).
For simplicity, the notation \(\bigO(\cdot)\) is used to conceal any constant that is independent of the dimension \(n\) and the iteration index \(k\).
For example, \(n\macheps + n^2\macheps^2\) can be expressed as \(\bigO(n\macheps)\) because \(n^2\macheps^2 \leq n\macheps\).

\subsection{Rounding error analysis for each step}
\label{subsec:analysis-step}
Taking rounding errors into account, we first present the results for each step of the \(k\)-th iteration of Algorithm~\ref{alg:pcg}, respectively.
In the following, \(\delta\) represents perturbations for individual values, while \(\Delta\) is employed to denote perturbations for vectors.
\begin{itemize}
\item Line~\ref{line:alpha} satisfies
\begin{align}
    \hatalphak = \frac{(\hatrk + \Delta\hatrk)\trans \hatqk}{\hatpk\trans A\hatpk+\delta e_{k}}(1+\delta r_k), \quad 
    \norm{\Delta\hatrk} 
    &\leq \bigO(n\macheps)\norm{\hatrk}, \label{eq:alphak}\\
    \abs{\delta r_k} 
    &\leq \bigO(\macheps),\label{eq:alphak-1} \\
    \abs{\delta e_k} &\leq \bigO(n\macheps)\norm{A}\norm{\hatpk}^2. \label{eq:alphak-2}
\end{align}
\item Line~\ref{line:x} satisfies
\begin{align}
    \hatxkone = \hatxk + \hatalphak\hatpk + \Dx_{k+1} + \Dbarx_{k+1}, \quad
    \norm{\Dx_{k+1}} &\leq \bigO(\macheps) \abs{\hatalphak}\norm{\hatpk}, \label{eq:xk} \\
    \norm{\Dbarx_{k+1}} &\leq \macheps \norm{\hatxk}. \label{eq:xk-1}
\end{align}
\item Line~\ref{line:r} satisfies
\begin{align}
    \hatrkone = \hatrk-\hatalphak A\hatpk + \Dr_{k+1}, \quad
    \norm{\Dr_{k+1}} &\leq \macheps\norm{\hatrk}
    +\bigO(n\macheps)\abs{\hatalphak}\norm{A}\norm{\hatpk}. \label{eq:hatrk1}
\end{align}
\item Lines~\ref{line:s} to~\ref{line:q} satisfy
\begin{align}
    \hatskone = M_L^{-1}\hatrkone + \Ds_{k+1}, \quad 
    \norm{\Ds_{k+1}}&\leq \epss\norm{M_L^{-1}\hatrkone}, \label{eq:hatsk1}\\
    \hatqkone = M_R^{-1}\hatskone + \Dq_{k+1}, \quad 
    \norm{\Dq_{k+1}} &\leq \epsq\norm{M_R^{-1}\hatskone}\label{eq:hatqk1},
\end{align}
where \(\epss\) and \(\epsq\) are parameters to bound the magnitude of errors, depending on how the left and right preconditioners are applied, respectively.
\item Line~\ref{line:beta} satisfies
\begin{align}
    \hatbetakone 
    = \frac{(\hatrkone+\Delta\hatrkone)\trans \hatqkone}{(\hatrk + \Delta\hatrk)\trans \hatqk}(1+\Dhatbeta_{k+1}),\quad
    \norm{\Delta\hatrkone} &\leq \bigO(n\macheps)\norm{\hatrkone}, \\
    \abs{\Dhatbeta_{k+1}} &\leq \bigO(\macheps). \label{eq:betak1}
\end{align}
\item Line~\ref{line:p} satisfies
\begin{align}
    \hatpkone = \hatqkone + \hatbetakone \hatpk + \Dp_{k+1}, \quad
    \norm{\Dp_{k+1}} &\leq \macheps\norm{\hatqkone}
    + \bigO(\macheps)\abs{\hatbetakone}\norm{\hatpk}. \label{eq:pk1}
\end{align}
\end{itemize}
Except for inequalities~\eqref{eq:hatsk1} and~\eqref{eq:hatqk1} that involve the application of preconditioners, these results follow straightforwardly from standard rounding error analysis under the assumption \(n\macheps<1\); see details in~\cite{H2002}.

As mentioned above, Greenbaum estimated the norm of the difference between the true residual \(\Tilde{r}_k \) and the recursively updated residual \(\hatrk\), and summarized the main result in~\cite[Theorem~2.2]{Greenbaum1997}.
Although her results were not originally formulated incorporating any preconditioner, it is straightforward to adapt them for the preconditioned case, in particular, for Algorithm~\ref{alg:pcg}. 
Here we present an equivalent formulation of this result for Algorithm~\ref{alg:pcg}, as well as the proof based on Greenbaum's ideas to emphasize that preconditioning does not affect the gap between \(\tilde r_k \) and \( \hatrk \).

\begin{theorem} \label{thm:Drk}
    Assume that \(\hatxk\) and \(\hatrk\) are generated by Algorithm~\ref{alg:pcg} satisfying~\eqref{eq:xk} and~\eqref{eq:hatrk1}.
    The difference between the true residual \(\Tilde{r}_k = b - A\hatxk\) and the recursively updated residual \(\hatrk\) satisfies
    \begin{equation*}
        \hatrk = \Tilde{r}_k + \Delta\Tilde{r}_k, \qquad
        \norm{\Delta\Tilde{r}_k}\leq \bigO(nk\macheps) \norm{A} \max_{j\leq k}\bigl(\norm{\hatx_j}, \norm{x}\bigr).
    \end{equation*}
\end{theorem}
\begin{proof}
    Consider the formulas~\eqref{eq:xk-1} and~\eqref{eq:hatrk1} to update the approximate solution and the recursively updated residual, respectively, and let \( \Delta \hatxk = \Delta x_k + \Delta \barxk\). Then 
    \begin{align*}
        \Delta \tilde r_k &= \hatrk - \tilde r_k\\
        &= \hat{r}_{k - 1} - \hat{\alpha}_{k - 1} A \hat{p}_{k - 1} + \Delta r_k 
        - b + A (\hat{x}_{k - 1} + \hat{\alpha}_{k - 1} \hat{p}_{k - 1} + \Delta \hatxk)\\
        &= \hat r_{k - 1} - (b - A \hat x_{k - 1}) + \Delta r_k + A \Delta \hat x_k.
    \end{align*}
    It is easy to see that, inductively, the residual gap \(\Delta \tilde r_k\) satisfies
    \begin{equation*}
        \Delta \tilde r_k = \Delta \tilde r_0 + \sum_{j = 1}^k A \Delta \hat x_j + \Delta r_j,
    \end{equation*}
    so that the quantities \(\hat \alpha_j\) and \( \hat p_j\), which use the preconditioned vectors in their formulas, cancel out for each \(j = 0, \dots, k - 1\).
    Taking the norm on both sides of the expression and subsequent applications of the triangle inequality yield 
    \begin{equation}\label{eq:thm:proof:resgap}
        \norm{\Delta \tilde r_k} 
        \leq \norm{\Delta \tilde r_0} + \sum_{j = 1}^k \norm{A} \norm{\Delta \hat x_j} + \norm{\Delta r_j}.
    \end{equation}
    In \cite[Lemma~2.1]{Greenbaum1997}, Greenbaum showed that the local errors in the sum satisfy
    \begin{equation}\label{eq:thm:proof:greenbaum}
        \sum_{j = 1}^k \norm{\Delta \hat x_j} \leq O(k \macheps) \max_{j \leq k} \norm{\hat x_j} \quad
        \text{ and } \quad
        \sum_{j = 1}^k \norm{\Delta r_j} \leq O(n k \macheps) \norm{A} \max_{j \leq k} \norm{\hat x_j}.
    \end{equation}
    Using standard rounding error analysis, it is not hard to show that 
    \begin{equation*}
        \norm{\Delta \tilde r_0} = \norm{\fl (b - A \hat x_0 ) - (b - A \hatx_0 )} \leq \norm{A} ( O(n \macheps)  \norm{\hatx_0} + \macheps \norm{x}). 
    \end{equation*}
    Substituting this and~\eqref{eq:thm:proof:greenbaum} into~\eqref{eq:thm:proof:resgap}, and taking the maximum between \( \norm{\hatx_j}, \norm{x} \) over all \(j = 1, \dots, k\) concludes the proof.
\end{proof}

This theorem presents the result that the recursively updated vector \(\hatrk\) is close to the true residual.
This means that the norm of the true residual can reach a small level for a sufficiently small norm of \(\hatrk\).

For convenience of the following analysis, we define the recurrence
\begin{equation} \label{eq:tildexk1}
    \bar x_{0} = x_0\quad\text{and}\quad \bar x_{k+1} = \bar x_{k} + \hatalphak\hatpk + \Dx_{k+1}.
\end{equation}
Combining this definition with~\eqref{eq:xk} allows us to rewrite \(\hatxkone\) as
\begin{equation} \label{eq:hatxk1-tildexk1}
    \hatxkone = \bar x_{k+1} + \sum_{j=0}^{k}\Dbarx_{j+1}, \qquad
    \normbig{\sum_{j=0}^{k}\Dbarx_{j+1}} \leq (k+1)\macheps\max_{j\leq k}\bigl(\norm{\hatx_j}\bigr).
\end{equation}
If we define \(\bar r_k = b - A \bar x_k\), then from Theorem~\ref{thm:Drk} and~\eqref{eq:hatxk1-tildexk1}, \(\hatrk\) can be written as
\begin{equation} \label{eq:hatrk-barrk}
    \hatrk = b - A\hatxk + \Delta\tilde r_k
    = b - A\bar x_k - A\sum_{j=0}^{k-1}\Dbarx_{j+1} + \Delta\tilde r_k
    = \bar r_k + \Delta\bar r_k,
\end{equation}
with \(\Delta\bar r_k := - A\sum_{j=0}^{k-1}\Dbarx_{j+1} + \Delta\tilde r_k\) satisfying
\begin{equation} \label{eq:normDbarr}
\begin{split}
    \norm{\Delta\bar r_k} &\leq \norm{A}\normbig{\sum_{j=0}^{k-1}\Dbarx_{j+1}} + \norm{\Delta\tilde r_k}
    \leq \bigO(nk\macheps)\norm{A}\max_{j\leq k}\bigl(\norm{\hatx_j},\norm{x}\bigr).
\end{split}
\end{equation}

Given \(M^{-1}=M_R^{-1}M_L^{-1}\), it is straightforward to see that in exact arithmetic, \(\alpha_k\) and \(\beta_{k + 1}\) can be equivalently expressed as
\begin{equation} \label{eq:exact-alpha}
\alpha_k = \frac{r_k\trans M^{-1}r_k}{p_k\trans Ap_k}\quad\text{and}\quad
\beta_{k + 1} = \frac{r_{k+1}\trans M^{-1}r_{k+1}}{r_k\trans M^{-1}r_k}.
\end{equation}
The following lemma provides a more explicit representation, which shows that the computed quantities \(\hatalphak\) and \(\hat{\beta}_{k + 1}\) satisfy perturbed versions of the above formulas.

\begin{lemma} \label{lem:eachstep-alphabeta}
    Assume that \(\hatrk\), \(\hatrkone\), \(\hatpk\), \(\hatalphak\), and \(\hatbetakone\) are computed by Algorithm~\ref{alg:pcg}.
    Define \(\epspre := \epss(1+\epsq) \kappa(M_R) \max_{y\neq 0} \bigl(\norm{M_L^{-1}y} \cdot\norm{M_R^{-1} y}/\normMinv{y}^2\bigr)+\bigl(\epsq+\bigO(n\macheps)\bigr) \kappa(M)^{1/2}\).
    If \(\bigO(n\macheps)\kappa(A)\leq 1/2\) and \(\epspre\leq 1/2\), then
    \begin{align}
        \hatalphak &= \frac{\hatrk\trans M^{-1} \hatrk}{\hatpk\trans A\hatpk}(1 + \Dalpha_k),
        &&\abs{\Dalpha_k} \leq \bigO(n\macheps)\kappa(A) + \bigO(\epspre), \label{eq:lem:normDalpha}\\
        \hatbetakone &= \frac{\hatrkone\trans M^{-1} \hatrkone}{\hatrk\trans M^{-1} \hatrk} (1+\Dbeta_{k+1}), 
        &&\abs{\Dbeta_{k+1}} \leq \bigO(\epspre) + \bigO(n\macheps). \label{eq:lem:betak1} 
    \end{align}
\end{lemma}

\begin{proof}
    We first prove~\eqref{eq:lem:normDalpha}.
    From~\eqref{eq:hatsk1} and~\eqref{eq:hatqk1}, \(\hatqkone\) satisfies 
    \begin{equation} \label{eq:lem-proof:hatqk1}
    \begin{split}
        \hatqkone = M_R^{-1}M_L^{-1}\hatrkone + \underbrace{M_R^{-1}\Ds_{k+1} + \Dq_{k+1}}_{=: \Delta\hatqkone},
    \end{split}
    \end{equation}
    where \(\Delta\hatqkone\) satisfies
    \begin{equation}
        \norm{\Delta\hatqkone}
        \leq\epss(1+\epsq)\norm{M_R^{-1}}\norm{M_L^{-1}\hatrkone}
        +\epsq\norm{M^{-1}\hatrkone}.
    \end{equation}
    From~\eqref{eq:alphak} and~\eqref{eq:lem-proof:hatqk1}, we obtain
    \begin{equation} \label{eq:lem:proof:hatalpha-0}
        \begin{split}
        \hatalphak &= \frac{\bigl(\hatrk + \Delta\hatrk\bigr)\trans \bigl(M_R^{-1}M_L^{-1}\hatrk+\Delta\hatqk\bigr)}{\hatpk\trans A\hatpk}
        \cdot\frac{1+\delta r_k}{1+\delta e_{k}/\hatpk\trans A\hatpk}.
        \end{split}
    \end{equation}
    On the right-hand side of~\eqref{eq:lem:proof:hatalpha-0}, from \(M^{-1}=M_R^{-1}M_L^{-1}\), \(\bigl(\hatrk + \Delta\hatrk\bigr)\trans \bigl(M_R^{-1}M_L^{-1}\hatrk+\Delta\hatqk\bigr)\) satisfies
    \begin{equation} \label{eq:lem:proof:hatzkhatsk}
        \begin{split}
            &\bigl(\hatrk + \Delta\hatrk\bigr)\trans \bigl(M_R^{-1}M_L^{-1}\hatrk+\Delta\hatqk\bigr) \\
            &= \hatrk\trans M^{-1}\hatrk + \hatrk\trans \Delta\hatqk + \Delta\hatrk\trans M^{-1}\hatrk + \Delta\hatrk\trans \Delta\hatqk \\
            &= \hatrk\trans M^{-1}\hatrk\biggl(1 + \underbrace{\frac{\hatrk\trans \Delta\hatqk + \Delta\hatrk\trans M^{-1}\hatrk + \Delta\hatrk\trans \Delta\hatqk}{\hatrk\trans M^{-1}\hatrk}}_{=:\Dalpha^{(1)}_k}\biggr).
        \end{split}
    \end{equation}
    This implies that~\eqref{eq:lem:proof:hatalpha-0} can be written as
    \begin{equation} \label{eq:lem:proof:hatalpha-1}
        \begin{split}
        \hatalphak 
        &= \frac{\hatrk\trans M^{-1}\hatrk}{\hatpk\trans A\hatpk}
        \underbrace{\biggl(\frac{1+\delta r_k}{1+\delta e_{k}/\hatpk\trans A\hatpk}
        + \Dalpha^{(1)}_k
        \cdot\frac{1+\delta r_k}{1+\delta e_{k}/\hatpk\trans A\hatpk}\biggr)}_{=:1+\Dalpha_k},
        \end{split}
    \end{equation}
    where \(\Dalpha_k\) satisfies, by the bound on \(\delta r_k\) in~\eqref{eq:alphak-1},
    \begin{equation} \label{eq:lem:proof:absDalphak-0}
    \begin{split}
        \abs{\Dalpha_k} &\leq \absbig{\frac{1+\delta r_k}{1+\delta e_{k}/\hatpk\trans A\hatpk}-1}
        + \abs{\Dalpha^{(1)}_k}\cdot \absbig{\frac{1+\delta r_k}{1+\delta e_{k}/\hatpk\trans A\hatpk}} \\
        &= \absbig{\frac{\delta r_k-\delta e_{k}/\hatpk\trans A\hatpk}{1+\delta e_{k}/\hatpk\trans A\hatpk}}
        + \abs{\Dalpha^{(1)}_k}\cdot \absbig{ 1+ \frac{1+\delta r_k}{1+\delta e_{k}/\hatpk\trans A\hatpk} - 1} \\
        &\leq \absbig{\frac{\delta r_k-\delta e_{k}/\hatpk\trans A\hatpk}{1+\delta e_{k}/\hatpk\trans A\hatpk}}
        + \abs{\Dalpha^{(1)}_k} \biggl(1+\absbig{\frac{\delta r_k-\delta e_{k}/\hatpk\trans A\hatpk}{1+\delta e_{k}/\hatpk\trans A\hatpk}}\biggr) \\
        &\leq \frac{\bigO(\macheps)+\abs{\delta e_{k}/\hatpk\trans A\hatpk}}{1-\abs{\delta e_{k}/\hatpk\trans A\hatpk}}
        + \abs{\Dalpha^{(1)}_k}\biggl(1+\frac{\bigO(\macheps)+\abs{\delta e_{k}/\hatpk\trans A\hatpk}}{1-\abs{\delta e_{k}/\hatpk\trans A\hatpk}}\biggr).
    \end{split}
    \end{equation}
    To estimate \(\abs{\Dalpha_k}\), we only need to bound \(\abs{\delta e_{k}/\hatpk\trans A\hatpk}\) and \(\abs{\Dalpha^{(1)}_k}\) involved in~\eqref{eq:lem:proof:absDalphak-0}.
    Recall that the Rayleigh quotient $\rho_A(y) := \langle A y,y \rangle / \norm{y}^2 $ of a symmetric positive definite matrix $A$ and a nonzero vector $y \in \mathbb R^n$ satisfies
    \[
        \lamdmin(A) = \min_{z \neq 0} \rho_A(z) \leq \rho_A(y) \leq \max_{z \neq 0} \rho_A(z) = \lamdmax(A).
    \]
    Hence, 
    \begin{equation} \label{eq:lem:proof:normAnormpk2/pkApk}
        \frac{\norm{A}\norm{\hatpk}^2}{\hatpk\trans A\hatpk}
        = \frac{\norm{A}}{\hatpk\trans A\hatpk/(\norm{\hatpk}^2)}
        \leq \frac{\norm{A}}{\lamdmin(A)}
        = \kappa(A),
    \end{equation}
    which implies that \(\abs{\delta e_{k}/\hatpk\trans A\hatpk}\leq \bigO(n\macheps)\kappa(A)\) by using the bound on \(\delta e_k\) in~\eqref{eq:alphak-2}.
    Applying the Cauchy--Schwarz inequality and using the bounds on \(\Delta\hatrk\), \(\Delta\hatqk\) in~\eqref{eq:alphak}, \eqref{eq:lem-proof:hatqk1}, yields the following bound for \(\Dalpha^{(1)}_k\), defined in~\eqref{eq:lem:proof:hatzkhatsk}:
    \begin{equation} \label{eq:lem:proof:Dalpha1k}
    \begin{split}
        \abs{\Dalpha^{(1)}_k}
        &\leq \frac{\epss(1+\epsq)\norm{M_R^{-1}}\norm{M_L^{-1}\hatrk}\norm{\hatrk}
        +(\epsq+\bigO(n\macheps))\norm{M^{-1}\hatrk}\norm{\hatrk}}{\hatrk\trans M^{-1}\hatrk} \\
        &\leq \epss(1+\epsq) \max_{y\neq 0} \biggl(\frac{\norm{M_L^{-1}y} \cdot\norm{M_R^{-1}}\cdot\norm{y}}{\normMinv{y}^2}\biggr) \\
        &\quad+\bigl(\epsq+\bigO(n\macheps)\bigr) \max_{y\neq 0} \biggl(\frac{\norm{M^{-1}y} \cdot\norm{y}}{\normMinv{y}^2}\biggr)\\
        &\leq \epss(1+\epsq) \max_{y\neq 0} \biggl(\frac{\norm{M_L^{-1}y} \cdot\norm{M_R^{-1}}\cdot\norm{M_R}\cdot\norm{M_R^{-1} y}}{\normMinv{y}^2}\biggr) \\
        &\quad+\bigl(\epsq+\bigO(n\macheps)\bigr) \max_{y\neq 0} \biggl(\frac{\norm{M^{-1/2}}\cdot\norm{M^{-1/2}y} \cdot\norm{M^{1/2}}\cdot\norm{M^{-1/2}y}}{\normMinv{y}^2}\biggr)\\
        &\leq \epss(1+\epsq) \kappa(M_R) \max_{y\neq 0} \biggl(\frac{\norm{M_L^{-1}y} \cdot\norm{M_R^{-1} y}}{\normMinv{y}^2}\biggr)
        +\bigl(\epsq+\bigO(n\macheps)\bigr) \kappa(M)^{1/2}\\
        &=\epspre,
    \end{split}
    \end{equation}
    where we eliminate the dependence on \(k\) by substituting \(\hatrk\) with the maximum taken over all possible vectors \(y\).
    Combining the bounds on \(\abs{\delta e_{k}/\hatpk\trans A\hatpk}\) and \(\abs{\Dalpha^{(1)}_k}\) with~\eqref{eq:lem:proof:absDalphak-0}, and using the assumptions \( \bigO(n\macheps) \kappa(A) \leq 1/2\) and \(\bigO(n\macheps) < 1\) gives the desired bound for \(\abs{\Dalpha_k}\), and thus proves~\eqref{eq:lem:normDalpha}.

    We will now prove~\eqref{eq:lem:betak1}.
    Substituting the expressions given in~\eqref{eq:lem-proof:hatqk1} for \(\hatqkone\) and \(\hatqk\) involved in~\eqref{eq:betak1},
    and then using~\eqref{eq:lem:proof:hatzkhatsk}, we have
    \begin{equation} \label{eq:lem:proof:betak1}
        \begin{split}
            \hatbetakone &= \frac{\bigl(\hatrkone + \Delta\hatrkone\bigr)\trans \bigl(M_R^{-1}M_L^{-1}\hatrkone+\Delta\hatqkone\bigr)}{\bigl(\hatrk + \Delta\hatrk\bigr)\trans \bigl(M_R^{-1}M_L^{-1}\hatrk+\Delta\hatqk\bigr)} (1+\Dhatbeta_{k+1}) \\
            &= \frac{\hatrkone\trans M^{-1} \hatrkone}{\hatrk\trans M^{-1} \hatrk} 
            \cdot\underbrace{\frac{1 + \Dalpha^{(1)}_{k+1}}{1 + \Dalpha^{(1)}_{k}}(1+\Dhatbeta_{k+1})}_{=:1+\Dbeta_{k+1}}.
        \end{split}
    \end{equation}
    Together with the bound on \(\abs{\Dalpha^{(1)}_k}\) in~\eqref{eq:lem:proof:Dalpha1k} and our assumptions on \( \epspre \), \(\Dbeta_{k+1}\) can be bounded as
    \begin{equation}
        \begin{split}
            \abs{\Dbeta_{k+1}} &= \absbig{\frac{1 + \Dalpha^{(1)}_{k+1}}{1 + \Dalpha^{(1)}_{k}}(1+\Dhatbeta_{k+1})-1} \\
            &= \absbig{\frac{\Dalpha^{(1)}_{k+1} + \Dhatbeta_{k+1} + \Dalpha^{(1)}_{k+1}\Dhatbeta_{k+1} - \Dalpha^{(1)}_{k}}{1 + \Dalpha^{(1)}_{k}}} \\
            &\leq \bigO(\epspre) + \bigO(n\macheps),
        \end{split}
    \end{equation}
    which proves~\eqref{eq:lem:betak1}.
\end{proof}

Lemma~\ref{lem:eachstep-alphabeta} provides a clearer expression for the computed results of \(\hatalphak\) and \(\hatbetakone\).
Now we consider \(\norm{p_{k+1}}\).
From the symmetric positive definite property of \(M^{-1}\) and \(M^{-1/2}\), we have \(M^{-1} = M^{-1/2}M^{-1/2}=(M^{-1/2})\trans M^{-1/2}\) and then
\begin{equation} \label{eq:lem-proof:bound-norm-Minvr0}
\begin{split}
    \norm{M^{-1}\hatr_0}
    &\leq \norm{M^{-1/2}}\norm{M^{-1/2}y} \\
    &= \sqrt{\lambda_{\max}(M^{-1})} \sqrt{y\trans M^{-1} y} \\
    &= \norm{M^{-1}}^{1/2}\normMinv{y}.
\end{split}
\end{equation}
In exact arithmetic, it is not hard to verify that
\begin{equation}
    p_{k+1} = M^{-1}r_{k+1}+\sum_{j=0}^k\frac{\normMinv{r_{k+1}}^2}{\normMinv{r_{j}}^2} M^{-1}r_{j},
\end{equation}
and thus, together with~\eqref{eq:lem-proof:bound-norm-Minvr0},
\begin{equation}
\begin{split}
    \norm{p_{k+1}} &\leq \norm{M^{-1}}^{1/2}\normMinv{r_{k+1}}
    + \sum_{j=0}^k\frac{\normMinv{r_{k+1}}^2}{\normMinv{r_{j}}^2} \norm{M^{-1}}^{1/2}\normMinv{r_{j}} \\
    &= \bigl(1+\sum_{j=0}^k\frac{\normMinv{r_{k+1}}}{\normMinv{r_{j}}}\bigr) \norm{M^{-1}}^{1/2}\normMinv{r_{k+1}}.
\end{split}
\end{equation}
Utilizing the formula for \(\hatbetakone\) as presented in Lemma~\ref{lem:eachstep-alphabeta}, we establish the relationship between the computed results \(\norm{\hatpkone}\) and \(\norm{\hatrkone}\) in the following lemma.

\begin{lemma} \label{lem:eachstep}
    Assume that \(\hatrkone\) and \(\hatpkone\) are computed by Algorithm~\ref{alg:pcg}, and \(\epspre\) is defined in Lemma~\ref{lem:eachstep-alphabeta}.
    If \(\bigO(n\macheps)\kappa(A)\leq 1/2\) and \(\epspre\leq 1/2\), then
    \begin{equation} \label{eq:lem:normpk-nokassump}
    \begin{split}
        \norm{\hatpkone}
        &\leq \biggl(1+\biggl(1+\bigO(\epspre)+ \bigO(n\macheps)\biggr)^{k+1}\biggl(\bigO\bigl((k+1)\epspre\bigr) + \bigO\bigl(n(k+1)\macheps\bigr)\biggr)\biggr)\\
        &\quad\cdot\biggl(1+\sum_{j=0}^{k}\frac{\normMinv{\hatr_{k+1}}}{\normMinv{\hatr_j}} \biggr) \norm{M^{-1}}^{\frac{1}{2}} \normMinv{\hatr_{k+1}}.
    \end{split}
    \end{equation}
    Furthermore, if \(\bigO\bigl((k+1)\epspre\bigr) + \bigO\bigl(n(k+1)\macheps\bigr)<1\), then
    \begin{align}
        \norm{\hatpkone}
        &\leq \biggl(1+\bigO\bigl((k+1)\epspre\bigr) + \bigO\bigl(n(k+1)\macheps\bigr)\biggr) \notag\\
        &\quad\cdot\biggl(1+\sum_{j=0}^{k}\frac{\normMinv{\hatr_{k+1}}}{\normMinv{\hatr_j}} \biggr) \norm{M^{-1}}^{\frac{1}{2}} \normMinv{\hatr_{k+1}}, \label{eq:lem:normpk-0} \\
        \macheps\norm{\hatpkone} &\leq \bigO(\macheps) \biggl(1+\sum_{j=0}^{k}\frac{\normMinv{\hatr_{k+1}}}{\normMinv{\hatr_j}}\biggr)
         \norm{M^{-1}}^{\frac{1}{2}} \normMinv{\hatr_{k+1}},   \label{eq:lem:u-normpk}
    \end{align}
    and
    \begin{equation} \label{eq:lem:Dpk1norm-0}
    \begin{split}
        \norm{\Dp_{k+1}} &\leq \bigO(\macheps) \biggl(1+\sum_{j=0}^{k}\frac{\normMinv{\hatr_{k+1}}}{\normMinv{\hatr_j}}\biggr)
         \norm{M^{-1}}^{\frac{1}{2}} \normMinv{\hatr_{k+1}}.
    \end{split}
    \end{equation}
\end{lemma}

\begin{proof}
    We will prove~\eqref{eq:lem:normpk-nokassump} by induction.
    Before proving the base case, using~\eqref{eq:lem-proof:bound-norm-Minvr0}, we first show that
    \begin{equation} \label{eq:lem:kappaMRbound}
    \begin{split}
        \frac{\norm{M_R^{-1}}\cdot\norm{M_L^{-1}\hatr_0}}{\norm{M^{-1}}^{\frac{1}{2}} \normMinv{\hatr_0}}
        &= \frac{\norm{M_R^{-1}}\cdot\norm{\hatr_0}\cdot\norm{M_L^{-1}\hatr_0}}{\norm{M^{-1}}^{\frac{1}{2}}\norm{\hatr_0}\cdot \normMinv{\hatr_0}}\\
        &\leq \frac{\norm{M_R^{-1}}\cdot\norm{M_R}\cdot\norm{M_R^{-1}\hatr_0}\cdot\norm{M_L^{-1}\hatr_0}}{\normMinv{\hatr_0}^2}\\
        &\leq \kappa(M_R) \max_{y}\biggl(\frac{\norm{M_R^{-1}y}\cdot\norm{M_L^{-1}y}}{\normMinv{y}^2}\biggr),
    \end{split}
    \end{equation}
    and
    \begin{equation}\label{eq:lem:kappaMbound}
    \begin{split}
        \frac{\norm{M^{-1}\hatr_0}}{\norm{M^{-1}}^{\frac{1}{2}} \normMinv{\hatr_0}}
        &= \frac{\norm{\hatr_0}\cdot\norm{M^{-1/2}}\cdot\norm{M^{-1/2}\hatr_0}}{\norm{M^{-1}}^{\frac{1}{2}}\norm{\hatr_0} \cdot\normMinv{\hatr_0}}
        \leq \kappa(M)^{1/2}.
    \end{split}
    \end{equation}
    For the base case, i.e., \(k+1=0\), due to~\eqref{eq:lem-proof:hatqk1}, \eqref{eq:lem-proof:bound-norm-Minvr0}, \eqref{eq:lem:kappaMRbound}, and~\eqref{eq:lem:kappaMbound},
    we have
    \begin{equation} \label{eq:lem:proof:normp0}
    \begin{split}
        \norm{\hatp_0}&=\norm{\hatq_0} \\
        &\leq \norm{M_R^{-1}M_L^{-1}\hatr_0} + \norm{\Delta\hat{q}_0} \\
        &\leq \norm{M^{-1}\hatr_0}
        +\epss(1+\epsq)\norm{M_R^{-1}}\norm{M_L^{-1}\hatr_0}
        +\epsq\norm{M^{-1}\hatr_0} \\
        &\leq \norm{M^{-1}}^{\frac{1}{2}} \normMinv{\hatr_0} 
        +\frac{\epss(1+\epsq)\norm{M_R^{-1}}\norm{M_L^{-1}\hatr_0}
        +\epsq\norm{M^{-1}\hatr_0}}{\norm{M^{-1}}^{\frac{1}{2}} \normMinv{\hatr_0}}\norm{M^{-1}}^{\frac{1}{2}} \normMinv{\hatr_0} \\
        &\leq \bigl(1+\epspreq\bigr) \norm{M^{-1}}^{\frac{1}{2}} \normMinv{\hatr_0},
    \end{split}
    \end{equation}
    where for the last inequality we have used~\eqref{eq:lem:kappaMRbound} and~\eqref{eq:lem:kappaMbound}.
    Then assuming that~\eqref{eq:lem:normpk-nokassump} holds for \(k+1=i\), then we aim to prove that it also holds for \(k+1=i+1\).
    From~\eqref{eq:pk1} and~\eqref{eq:lem:betak1}, we obtain
    \begin{equation} \label{eq:lem:proof:normpi1-0}
    \begin{split}
        \norm{\hatp_{i+1}} &\leq (1+\macheps)\norm{\hatq_{i+1}} + \bigl(1+\bigO(\macheps)\bigr)\norm{\hatbeta_{i+1} \hatp_i} \\
        &\leq (1+\macheps)\norm{\hatq_{i+1}} + \biggl(1+\bigO(\epspre) + \bigO(n\macheps)\biggr)\frac{\hatr_{i+1}\trans M^{-1}\hatr_{i+1}}{\hatr_i\trans M^{-1}\hatr_i} \norm{\hatp_i}.
    \end{split}
    \end{equation}
    Similarly to~\eqref{eq:lem:proof:normp0}, the first term involved in the right-hand side of~\eqref{eq:lem:proof:normpi1-0} satisfies
    \begin{equation} \label{eq:lem:proof:righthand-normpi1-1}
        (1+\macheps)\norm{\hatq_{i+1}}
        \leq \bigl(1+\bigO(\macheps) + \epspreq\bigr) \norm{M^{-1}}^{\frac{1}{2}}\normMinv{\hatr_{i+1}}.
    \end{equation}
    By using the induction hypothesis on \(\norm{\hatp_i}\), the second term involved in the right-hand side of~\eqref{eq:lem:proof:normpi1-0} satisfies
    \begin{equation} \label{eq:lem:proof:righthand-normpi1-2-0}
        \begin{split}
            &\biggl(1+\bigO(\epspre)+ \bigO(n\macheps)\biggr)\frac{\hatr_{i+1}\trans M^{-1}\hatr_{i+1}}{\hatr_i\trans M^{-1}\hatr_i} \norm{\hatp_i} \\
            &\leq \biggl(1+\bigO(\epspre)+ \bigO(n\macheps)\biggr)
            \biggl(1+\biggl(1+\bigO(\epspre)+ \bigO(n\macheps)\biggr)^i\bigl(\bigO(i\epspre) + \bigO\bigl(ni\macheps\bigr)\bigr)\biggr) \\
            &\quad\cdot
            \biggl(1+\sum_{j=0}^{i-1}\frac{\normMinv{\hatr_i}}{\normMinv{\hatr_j}}\biggr)\frac{\normMinv{\hatr_{i+1}}^2}{\normMinv{\hatr_i}} \norm{M^{-1}}^{\frac{1}{2}}.
        \end{split}
    \end{equation}
    In the right-hand side of~\eqref{eq:lem:proof:righthand-normpi1-2-0}, the first two terms can be bounded by
    \begin{equation} \label{eq:lem:proof:righthand-normpi1-2-1}
    \begin{split}
        &\biggl(1+\bigO(\epspre)+ \bigO(n\macheps)\biggr)
            \biggl(1+\biggl(1+\bigO(\epspre)+ \bigO(n\macheps)\biggr)^i\bigl(\bigO(i\epspre) + \bigO\bigl(ni\macheps)\bigr)\biggr) \\
        &= 1+\bigO(\epspre)+ \bigO(n\macheps) + \biggl(1+\bigO(\epspre)+ \bigO(n\macheps)\biggr)^{i+1}\biggl(\bigO\bigl(i\epspre\bigr) + \bigO\bigl(ni\macheps\bigr)\biggr) \\
        &\leq 1+\biggl(1+\bigO(\epspre)+ \bigO(n\macheps)\biggr)^{i+1}\biggl(\bigO\bigl((i+1)\epspre\bigr) + \bigO\bigl(n(i+1)\macheps\bigr)\biggr).
    \end{split}
    \end{equation}
    The remaining terms can be expressed as
    \begin{equation} \label{eq:lem:proof:righthand-normpi1-2-2}
    \begin{split}
        \biggl(1+\sum_{j=0}^{i-1}\frac{\normMinv{\hatr_i}}{\normMinv{\hatr_j}}\biggr)\frac{\normMinv{\hatr_{i+1}}^2}{\normMinv{\hatr_i}} \norm{M^{-1}}^{\frac{1}{2}}
        &= \sum_{j=0}^{i}\frac{\normMinv{\hatr_i}}{\normMinv{\hatr_j}}\cdot\frac{\normMinv{\hatr_{i+1}}^2}{\normMinv{\hatr_i}} \norm{M^{-1}}^{\frac{1}{2}} \\
        &= \sum_{j=0}^{i}\frac{\normMinv{\hatr_{i+1}}}{\normMinv{\hatr_j}}\cdot\norm{M^{-1}}^{\frac{1}{2}} \normMinv{\hatr_{i+1}}.
    \end{split}
    \end{equation}
    Substituting these two terms involved in the right-hand side of~\eqref{eq:lem:proof:righthand-normpi1-2-0} by~\eqref{eq:lem:proof:righthand-normpi1-2-1} and~\eqref{eq:lem:proof:righthand-normpi1-2-2}, we have
    \begin{equation} \label{eq:lem:proof:righthand-normpi1-2}
        \begin{split}
            &\biggl(1+\bigO(\epspre)+ \bigO(n\macheps)\biggr)\frac{\hatr_{i+1}\trans M^{-1}\hatr_{i+1}}{\hatr_i\trans M^{-1}\hatr_i} \norm{\hatp_i} \\
            &\leq \biggl(1+\biggl(1+\bigO(\epspre)+ \bigO(n\macheps)\biggr)^{i+1}\biggl(\bigO\bigl((i+1)\epspre\bigr) + \bigO\bigl(n(i+1)\macheps\bigr)\biggr)\biggr) \\
            &\quad\cdot
            \sum_{j=0}^{i}\frac{\normMinv{\hatr_{i+1}}}{\normMinv{\hatr_j}}\cdot\norm{M^{-1}}^{\frac{1}{2}} \normMinv{\hatr_{i+1}}.
        \end{split}
    \end{equation}
    Substituting the terms involved in the right-hand side of~\eqref{eq:lem:proof:normpi1-0} by~\eqref{eq:lem:proof:righthand-normpi1-1} and~\eqref{eq:lem:proof:righthand-normpi1-2}, we derive
    \begin{equation} \label{eq:lem:proof:normpi1}
    \begin{split}
        \norm{\hatp_{i+1}} 
        &\leq \biggl(1+\biggl(1+\bigO(\epspre)+ \bigO(n\macheps)\biggr)^{i+1}\biggl(\bigO\bigl((i+1)\epspre\bigr) + \bigO\bigl(n(i+1)\macheps\bigr)\biggr)\biggr)\\
        &\quad\cdot\biggl(1+\sum_{j=0}^{i}\frac{\normMinv{\hatr_{i+1}}}{\normMinv{\hatr_j}} \biggr) \norm{M^{-1}}^{\frac{1}{2}} \normMinv{\hatr_{i+1}}.
    \end{split}
    \end{equation}
    Hence, the induction is proved and~\eqref{eq:lem:normpk-nokassump} holds.

    With the assumption \(\bigO\bigl((k+1)\epspre\bigr) + \bigO\bigl(n(k+1)\macheps\bigr)<1\), we have
    \begin{equation}
        \biggl(1+\bigO(\epspre)+ \bigO(n\macheps)\biggr)^{k+1}
        \leq \exp\biggl\{(k+1)\bigl(\bigO(\epspre)+ \bigO(n\macheps)\bigr)\biggr\}
        < 3,    
    \end{equation}
    Together with~\eqref{eq:lem:normpk-nokassump}, this proves~\eqref{eq:lem:normpk-0}.
    Multiplying the right-hand side of~\eqref{eq:lem:normpk-0} by $\macheps$, and using the assumption again yields~\eqref{eq:lem:u-normpk} .
    Furthermore, comparing~\eqref{eq:pk1} and~\eqref{eq:lem:proof:normpi1-0}, it is easy to prove~\eqref{eq:lem:Dpk1norm-0}.
\end{proof}

\subsection{Local orthogonality}
\label{subsec:localorth}
In this subsection, we analyze the local orthogonality, i.e., \(\abs{\hatrkone\trans \hatpk}\).
It should be noted that \(r_{k+1}\trans p_k = 0\) and \(r_k\trans p_k = \normMinv{r_k}^2\) hold in exact arithmetic.

\begin{lemma} \label{lem:rkpk-rk1pk}
    Assume that \(\hatrk\) and \(\hatpk\) are generated by Algorithm~\ref{alg:pcg}, and let \(\epspre\) be as defined in Lemma~\ref{lem:eachstep-alphabeta}.
    If \(\bigO(n\macheps)\kappa(A)\leq 1/2\), \(\epspreq\leq 1/2\), and \(\bigO(k\,\epspre) + \bigO(nk\macheps)<1\), then
    \begin{equation} \label{eq:lem:rk1pk}
    \begin{split}
        \abs{\hatrkone\trans \hatpk}&\leq \exp\bigl\{k\bigO(n\macheps) + k\bigO\bigl(\epspre\bigr)\bigr\}
        \biggl(\bigO\bigl(nk\macheps\bigr)\kappa(A)+ \bigO\bigl(k\epspre\bigr) \\
        &\quad+ \bigO\bigl(\macheps\bigr)\kappa(M)^{\frac{1}{2}} 
        \sum_{l=0}^{k}\sum_{j=0}^{l}\frac{\normMinv{\hatr_{l}}}{\normMinv{\hatr_j}}\biggr)\normMinv{\hatr_{k}}^2 \\
        &\leq \biggl(\bigO\bigl(nk\macheps\bigr)\kappa(A)+ \bigO\bigl(k\epspre\bigr)\biggr)\normMinv{\hatr_{k}}^2 \\
        &\quad+ \bigO\bigl(\macheps\bigr)\kappa(M)^{\frac{1}{2}} 
        \biggl(k+1+\sum_{l=0}^{k}\sum_{j=0}^{l-1}\frac{\normMinv{\hatr_{l-1}}}{\normMinv{\hatr_j}}\biggr)\normMinv{\hatr_{k}}^2.
    \end{split}
    \end{equation}
\end{lemma}

\begin{proof}
    We will first give an expression for \(\hatrkone\trans\hatpk\).
    From~\eqref{eq:hatrk1}, \(\hatrkone\trans\hatpk\) can be written as
    \begin{equation} \label{eq:lem:proof:rk1pk-0}
    \begin{split}
        \hatrkone\trans\hatpk &= \hatrk\trans\hatpk-\hatalphak \hatpk\trans A\hatpk + \Dr_{k+1}\trans\hatpk \\
        &= \underbrace{\bigl(\hatrk\trans\hatpk-\normMinv{\hatrk}^2\bigr)}_{=: \Drp_{k}}
        + \bigl(\normMinv{\hatrk}^2-\hatalphak \hatpk\trans A\hatpk\bigr)
        + \Dr_{k+1}\trans\hatpk.
    \end{split}
    \end{equation}
    We then treat \(\normMinv{\hatrk}^2-\hatalphak \hatpk\trans A\hatpk\), \(\Dr_{k+1}\trans\hatpk\), and \(\Drp_{k}\) involved in the bound.
    By~\eqref{eq:lem:normDalpha}, we have
    \begin{equation} \label{eq:lem:proof:rk2-alphapkApk}
    \begin{split}
        \absbig{\normMinv{\hatrk}^2-\hatalphak\hatpk\trans A\hatpk}
        &= \absbig{\normMinv{\hatrk}^2 - \normMinv{\hatrk}^2 \cdot(1+\Dalpha_k)} \\
        &\leq \biggl(\bigO(n\macheps)\kappa(A) + \bigO(\epspre)\biggr)\normMinv{\hatrk}^2.
    \end{split}
    \end{equation}
    Using the bound on \(\Dr_{k+1}\) in~\eqref{eq:hatrk1}, \(\Dr_{k+1}\trans \hatpk\) can be bounded as
    \begin{equation} \label{eq:lem:proof:Drk1pk}
    \begin{split}
        \abs{\Dr_{k+1}\trans \hatpk}
        &\leq \macheps\norm{\hatrk}\norm{\hatpk}+\bigO(n\macheps)\abs{\hatalphak}\norm{A}\norm{\hatpk}^2 \\
        &\leq \macheps \bigl(\normMinv{\hatrk} \norm{M}^{\frac{1}{2}} \bigr) \norm{\hatpk} 
        +\bigO(n\macheps)\kappa(A)\normMinv{\hatrk}^2 \\
        &\leq \bigO(\macheps) \kappa(M)^{\frac{1}{2}}  \biggl(1+\sum_{j=0}^{k-1}\frac{\normMinv{\hatrk}}{\normMinv{\hatr_j}}\biggr)\normMinv{\hatrk}^2+\bigO(n\macheps)\kappa(A)\normMinv{\hatrk}^2,
    \end{split}
    \end{equation}
    where the second inequality is derived from the expression for \(\hatalphak\) in~\eqref{eq:lem:normDalpha},~\eqref{eq:lem:proof:normAnormpk2/pkApk}, and the assumptions above. The last inequality is obtained due to the bound on \(\macheps\norm{\hatpk}\) in~\eqref{eq:lem:u-normpk}.
    Using the expressions~\eqref{eq:pk1}, \eqref{eq:lem:betak1}, and~\eqref{eq:lem-proof:hatqk1} for \(\hatp_{k+1}\), \(\hatbeta_{k+1}\), and \(\hatqkone\), respectively, we notice that \(\Drp_{k}\) can be written as
    \begin{equation} \label{eq:lem:proof:ri1pi1-0}
        \begin{split}
            \Drp_{k}
            &= \hatr_{k}\trans \hatp_{k} - \normMinv{\hatr_{k}}^2 \\
            &= \hatr_{k}\trans (M^{-1}\hatr_{k}+\Delta\hatqk+\hatbeta_{k}\hatp_{k-1}+\Dp_{k}) - \normMinv{\hatr_{k}}^2 \\
            &= \hatr_{k}\trans \Delta\hatqk + \frac{\hatr_{k}\trans M^{-1}\hatr_{k}}{\hatr_{k-1}\trans M^{-1}\hatr_{k-1}}\bigl(1 + 
            \Dbeta_{k}\bigr)\hatr_{k}\trans\hatp_{k-1}
            + \hatr_{k}\trans\Dp_{k}.
        \end{split}
    \end{equation}
    On the right-hand side of~\eqref{eq:lem:proof:ri1pi1-0}, from the bound on \(\Delta\hatqkone\) in~\eqref{eq:lem-proof:hatqk1}, similarly to~\eqref{eq:lem:proof:Dalpha1k}, we have
    \begin{equation} \label{eq:lem:proof:ri1pi1-subproof-0}
    \begin{split}
        \abs{\hatr_{k}\trans \Delta\hatqk}
        &\leq \frac{\epss(1+\epsq)\norm{M_R^{-1}}\cdot\norm{M_L^{-1}\hatrk}\cdot\norm{\hatrk}
        +\epsq\norm{M^{-1}\hatrk}\cdot\norm{\hatrk}}{\normMinv{\hatrk}^2}\cdot \normMinv{\hatrk}^2 \\
        &\leq \epspreq\normMinv{\hatrk}^2,
    \end{split}
    \end{equation}
    and from the upper bound on \(\Dp_{k}\) shown in~\eqref{eq:lem:Dpk1norm-0},
    \begin{equation} \label{eq:lem:proof:rDp}
        \begin{split}
            \abs{\hatr_{k}\trans\Dp_{k}}&\leq \norm{\hatr_{k}} \cdot\norm{\Dp_{k}} \\
            &\leq \bigl(\normMinv{\hatr_{k}}\norm{M}^{\frac{1}{2}} \bigr) \norm{\Dp_{k}} \\
            &\leq \bigO\bigl(\macheps\bigr)\kappa(M)^{\frac{1}{2}}
            \biggl(1+\sum_{j=0}^{k-1}\frac{\normMinv{\hatr_{k}}}{\normMinv{\hatr_j}}\biggr) \normMinv{\hatr_{k}}^2 \\
            &= \bigO\bigl(\macheps\bigr)\kappa(M)^{\frac{1}{2}}
            \sum_{j=0}^{k}\frac{\normMinv{\hatr_{k}}}{\normMinv{\hatr_j}} \normMinv{\hatr_{k}}^2.
        \end{split}
    \end{equation}
    Substituting \(\hatr_{k}\trans \Delta\hatqk\) and \(\hatr_{k}\trans\Dp_{k}\) involved in~\eqref{eq:lem:proof:ri1pi1-0} by~\eqref{eq:lem:proof:ri1pi1-subproof-0} and~\eqref{eq:lem:proof:rDp}, \(\Drp_{k}\) can be further bounded by
    \begin{equation} \label{eq:lem:proof:ri1pi1-1}
        \begin{split}
            \abs{\Drp_{k}}
            &\leq \epspreq\normMinv{\hatr_{k}}^2 + \absbig{\frac{\normMinv{\hatr_{k}}^2}{\normMinv{\hatr_{k-1}}^2}\bigl(1 + \Dbeta_{k}\bigr)\hatr_{k}\trans\hatp_{k-1}} \\
            &\quad+ \bigO\bigl(\macheps\bigr)\kappa(M)^{\frac{1}{2}}
            \sum_{j=0}^{k}\frac{\normMinv{\hatr_{k}}}{\normMinv{\hatr_j}} \normMinv{\hatr_{k}}^2 \\
            &\leq \epspreq\normMinv{\hatr_{k}}^2 + \bigl(1 + \bigO\bigl(\epspre\bigr) + \bigO(n\macheps)\bigr)\frac{\normMinv{\hatr_{k}}^2}{\normMinv{\hatr_{k-1}}^2}\abs{\hatr_{k}\trans\hatp_{k-1}} \\
            &\quad+ \bigO\bigl(\macheps\bigr)\kappa(M)^{\frac{1}{2}}
            \sum_{j=0}^{k}\frac{\normMinv{\hatr_{k}}}{\normMinv{\hatr_j}} \normMinv{\hatr_{k}}^2.
        \end{split}
    \end{equation}
    
    Then combining~\eqref{eq:lem:proof:rk1pk-0} with the bound on \(\abs{\normMinv{\hatrk}^2-\hatalphak\hatpk\trans A\hatpk}\) in~\eqref{eq:lem:proof:rk2-alphapkApk}, the bound on \(\abs{\Dr_{k+1}\trans \hatpk}\) in~\eqref{eq:lem:proof:Drk1pk}, and the bound on \(\Drp_{k}\) in~\eqref{eq:lem:proof:ri1pi1-1}, \(\hatrkone\trans\hatpk\) can be bounded as
    \begin{equation} \label{eq:lem:proof:rk1pk}
    \begin{split}
        \abs{\hatrkone\trans\hatpk}
        &\leq \bigl(1 + \bigO\bigl(\epspre\bigr) + \bigO(n\macheps)\bigr)\frac{\normMinv{\hatr_{k}}^2}{\normMinv{\hatr_{k-1}}^2}\abs{\hatr_{k}\trans\hatp_{k-1}} \\
        &\quad+ \biggl(\bigO(n\macheps)\kappa(A) + \bigO\bigl(\epspre\bigr)\biggr)\normMinv{\hatrk}^2 \\
        &\quad+ \bigO(\macheps) \kappa(M)^{\frac{1}{2}} \biggl(1+\sum_{j=0}^{k-1}\frac{\normMinv{\hatrk}}{\normMinv{\hatr_j}}\biggr)\normMinv{\hatrk}^2.
    \end{split}
    \end{equation}

    From~\eqref{eq:lem:proof:rk1pk}, we find that the bound on \(\abs{\hatrkone\trans\hatpk}\) is determined by the bound on \(\abs{\hatr_{k}\trans\hatp_{k-1}}\).
    Thus, we will now prove the conclusion~\eqref{eq:lem:rk1pk} by induction.
    For the base case, i.e., \(k=0\), noticing \(\hatp_0=\hatq_0\), similarly to~\eqref{eq:lem:proof:Dalpha1k} derived by using~\eqref{eq:lem-proof:hatqk1}, we have
    \begin{equation} \label{eq:lem:proof:r0p0-r0Minv}
    \begin{split}
        \abs{\hatr_0\trans &\hatp_0 - \normMinv{\hatr_0}^2} \\
        &\leq \norm{\hatr_0}\cdot\norm{\Delta\hatq_0} \\
        &\leq \frac{\epss(1+\epsq)\norm{M_R^{-1}}\cdot\norm{M_L^{-1}\hatr_0}\cdot\norm{\hatr_0}
        +\epsq\norm{M^{-1}\hatr_0}\cdot\norm{\hatr_0}}{\normMinv{\hatr_0}^2}\cdot \normMinv{\hatr_0}^2\\
        &\leq \epspreq \normMinv{\hatr_0}^2,
    \end{split}
    \end{equation}
    where the second inequality is derived from \(\norm{\hatr_0}\leq\norm{M_R^{-1}\hatr_0}\cdot\norm{M_R}\).
    Together with~\eqref{eq:lem:proof:rk1pk-0}, \eqref{eq:lem:proof:rk2-alphapkApk}, and~\eqref{eq:lem:proof:Drk1pk}, then \(\hatr_1\trans\hatp_0\) can be written as
    \begin{equation*}
    \begin{split}
        \abs{\hatr_1\trans\hatp_0} &\leq \biggl(\bigO(n\macheps)\kappa(A) + \bigO\bigl(\epspre\bigr)\biggr)\normMinv{\hatr_0}^2
        + \bigO(\macheps) \kappa(M)^{\frac{1}{2}} \normMinv{\hatr_0}^2,
    \end{split}
    \end{equation*}
    which gives the base case.
    
    Then we assume that~\eqref{eq:lem:rk1pk} holds for \(k=i\), and aim to prove that it also holds for \(k=i+1\).
    Using~\eqref{eq:lem:betak1} and the induction hypothesis on \(\hatr_{i+1}\trans\hatp_{i}\), the first term involved in the right-hand side of~\eqref{eq:lem:proof:rk1pk} for \(k=i+1\) can be bounded by
    \begin{equation}  \label{eq:lem:proof:ri1pi1-subproof-1}
    \begin{split}
        &\bigl(1 + \bigO(n\macheps) + \bigO\bigl(\epspre\bigr)\bigr)\frac{\normMinv{\hatr_{i+1}}^2}{\normMinv{\hatr_{i}}^2} \abs{\hatr_{i+1}\trans\hatp_{i}} \\
        &\leq \bigl(1 + \bigO(n\macheps) + \bigO\bigl(\epspre\bigr)\bigr)
        \exp\bigl\{i\,\bigO(n\macheps) + i\,\bigO\bigl(\epspre\bigr)\bigr\}
        \biggl(\bigO\bigl(ni\macheps\bigr)\kappa(A)+ \bigO\bigl(i\epspre\bigr) \\
        &\quad+ \bigO\bigl(\macheps\bigr)\kappa(M)^{\frac{1}{2}} 
        \sum_{l=0}^i\sum_{j=0}^{l}\frac{\normMinv{\hatr_{l}}}{\normMinv{\hatr_j}}\biggr)\normMinv{\hatr_{i+1}}^2 \\
        &\leq \exp\bigl\{(i+1)\bigO(n\macheps) + (i+1)\bigO\bigl(\epspre\bigr)\bigr\}
        \biggl(\bigO\bigl(ni\macheps\bigr)\kappa(A)+ \bigO\bigl(i\epspre\bigr) \\
        &\quad+ \bigO\bigl(\macheps\bigr)\kappa(M)^{\frac{1}{2}} 
        \sum_{l=0}^i\sum_{j=0}^{l}\frac{\normMinv{\hatr_{l}}}{\normMinv{\hatr_j}}\biggr)\normMinv{\hatr_{i+1}}^2,
    \end{split}
    \end{equation}
    where the last inequality uses the inequality \(1 + \bigO(n\macheps) + \bigO\bigl(\epspre\bigr)\leq \exp\bigl\{\bigO(n\macheps) + \bigO\bigl(\epspre\bigr)\bigr\}\).
    Together with~\eqref{eq:lem:proof:rk1pk} for \(k = i+1\), and observing that \(\exp\bigl\{(i+1)\bigO(n\macheps) + (i+1)\bigO\bigl(\epspre\bigr)\bigr\}\geq 1\), we obtain
    \begin{equation}
    \begin{split}
        \abs{\hatr_{i+2}\trans\hatp_{i+1}}
        &\leq \exp\bigl\{(i+1)\bigO(n\macheps) + (i+1)\bigO\bigl(\epspre\bigr)\bigr\}
        \biggl(\bigO\bigl(n(i+1)\macheps\bigr)\kappa(A)+ \bigO\bigl((i+1)\epspre\bigr) \\
        &\quad+ \bigO\bigl(\macheps\bigr)\kappa(M)^{\frac{1}{2}} 
        \sum_{l=0}^{i+1}\sum_{j=0}^{l}\frac{\normMinv{\hatr_{l}}}{\normMinv{\hatr_j}}\biggr)\normMinv{\hatr_{i+1}}^2.
    \end{split}
    \end{equation}
    Then we can conclude the proof by induction on \(k\), where \(\exp\bigl\{k\bigO(n\macheps) + k\bigO\bigl(\epspre\bigr)\bigr\}\leq 3\) by the assumption \(\bigO(k\,\epspre) + \bigO(nk\macheps)<1\).
\end{proof}

From Lemma~\ref{lem:rkpk-rk1pk}, the local orthogonality, i.e., \(\abs{\hatrkone\trans \hatpk}\), depends on \(\sum_{j=0}^{l-1}\bigl(\normMinv{\hatr_l}/\normMinv{\hatr_j}\bigr)\) and \(\normMinv{\hatrk}\), where \(l\leq k\).
Notice that, from~\cite[Equation (4.34)]{ST2005},
\begin{equation*}
\begin{split}
    \bigO\bigl(\macheps\bigr)\sum_{j=0}^{l-1}\frac{\normMinv{\hatr_l}}{\normMinv{\hatr_j}}
    &\leq \bigO\bigl(\macheps\bigr)\kappa(M)^{\frac{1}{2}} 
    \sum_{j=0}^{l-1}\frac{\norm{\hatr_l}}{\norm{\hatr_j}}
    \leq \bigO\bigl(k\macheps\bigr)\kappa(M)^{\frac{1}{2}} \kappa(A)^{\frac{1}{2}}.
\end{split}
\end{equation*}
Applying the above estimate to \(\sum_{j=0}^{l-1}\bigl(\normMinv{\hatr_l}/\normMinv{\hatr_j}\bigr)\) yields a bound for \(\abs{\hatrkone\trans \hatpk}\) from Lemma~\ref{lem:rkpk-rk1pk} that is consistent with the one presented in~\cite[Theorem 4.3]{ST2005}.
As we will show in Section~\ref{subsec:back-forward-error}, when the gap between \(f(\barxkone)\) and \(f(\barxk)\) becomes small enough, one expects \(\normMinv{\hatrk}/\normMinv{\hatr_j}\) will typically remain \(\bigO(1)\) throughout the iterations, because \(\normMinv{\hatr_j}\rightarrow \normMinv{\hatrk}\) as \(j\) approaches \(k\).
Therefore, we expect the term \(\sum_{j=0}^{l-1}\bigl(\normMinv{\hatr_l}/\normMinv{\hatr_j}\bigr)\) to be bounded by \(\bigO(k)\) rather than by \(\bigO(k) \kappa(M)^{\frac{1}{2}}\kappa(A)^{\frac{1}{2}}\) for some specific \(k\).
Consequently, when bounding \( \abs{\hatrkone\trans \hatpk}\) it is preferable to retain the term \(\sum_{j=0}^{l-1}\bigl(\normMinv{\hatr_l}/\normMinv{\hatr_j}\bigr)\). 

\subsection{Backward and forward error analysis of PCG}
\label{subsec:back-forward-error}
We are now prepared to present the backward and forward error analysis of CG.
As mentioned in Section~\ref{sec:introduction}, the CG algorithm can be regarded as an optimization method used to minimize the quadratic function \eqref{eq:quadratic-function}.
In both mathematical and numerical contexts, the difference \(f(\barxk) - f(\barxkone)\) can become sufficiently small in a specific iteration; see the detailed proof in Lemma~\ref{lem:Dfsmall}.
In the following, we will show that the backward and forward errors can achieve an ideal level in this particular iteration.
The proof involves four steps:
\begin{enumerate}
    \item First, we will establish the relationship between \(f(\barxk) - f(\barxkone)\) and \(\normMinv{\hatrk}\) in Lemma~\ref{lem:f-normrk}.\\
    \item Subsequently, in Lemma~\ref{lem:Dfsmall}, we will present conditions for the existence of an iteration in which \(f(\barxk) - f(\barxkone)\) can be sufficiently close.\\
    \item Based on Lemmas~\ref{lem:f-normrk} and~\ref{lem:Dfsmall}, we will give an initial estimate of the size of \(\hatrk\) when \(f(\barxk) - f(\barxkone)\) becomes small enough in Lemma~\ref{lem:backward-forward-err}.\\
    \item Finally, we will present our main results in Theorem~\ref{thm:backward-forward-err}, which state that the backward and forward errors of CG can also be very small when \(f(\barxk) - f(\barxkone)\) is sufficiently small.
\end{enumerate}

Before presenting Lemma~\ref{lem:f-normrk}, we will establish the fact that, in exact arithmetic, \(f(\barxkone) - f(\barxk)\) can be expressed as a quadratic function of \(\normMinv{\hatrk}\).
Using the definition of \(f\) in~\eqref{eq:quadratic-function} together with \(x_{k+1} = x_k + \alpha_k p_k\), and employing the identity \(r_k\trans p_k = \normMinv{r_k}^2\) established before Lemma~\ref{lem:rkpk-rk1pk} and~\eqref{eq:exact-alpha}, we obtain
\begin{equation} \label{eq:exact-arithmetic-difference}
\begin{split}
    f(x_{k+1}) - f(x_k)
    &= \frac{1}{2}x_{k+1}\trans A x_{k+1} - x_{k+1}\trans b
       - \frac{1}{2}x_k\trans A x_k + x_k\trans b \\
    &= \frac{1}{2}\alpha_k^2 p_k\trans A p_k + \alpha_k p_k\trans A x_k
       - \alpha_k p_k\trans b \\
    &= \frac{1}{2}\alpha_k^2 p_k\trans A p_k - \alpha_k p_k\trans r_k \\
    &= -\frac{1}{2}\alpha_k \normMinv{r_k}^2.
\end{split}
\end{equation}
Lemma 4 then establishes a finite precision analogue, namely, that  \(f(\barxkone) - f(\barxk)\) satisfies a perturbed version of ~\eqref{eq:exact-arithmetic-difference}.

\begin{lemma} \label{lem:f-normrk}
    Assume that \(\barxk\) and \(\barxkone\) satisfy~\eqref{eq:tildexk1}, and \(\hatalphak\) satisfies~\eqref{eq:alphak}.
    Let \(\epspre\) be as defined in Lemma~\ref{lem:eachstep-alphabeta}.
    If \(\bigO(n\macheps)\kappa(A)\leq 1/2\) and \(\bigO\bigl(k\epspre\bigr) + \bigO(nk\macheps)\leq 1/2\), then
    \begin{equation}
    \begin{split}
        f(\barxkone)-f(\barxk) 
        &= - \frac{1}{2}\hatalphak(1+\Dhatalpha_k)\normMinv{\hatrk}^2 + \delta f_k,
    \end{split}
    \end{equation}
    where \(\Dhatalpha_k\) and \(\delta f_k\) satisfy, respectively,
    \begin{align}
        \abs{\Dhatalpha_k} &\leq \bigO\bigl(nk\macheps\bigr)\kappa(A)+ \bigO\bigl(k\epspre\bigr)
        + \bigO\bigl(\macheps\bigr)\kappa(M)^{\frac{1}{2}} 
        \biggl(k+1+\sum_{l=0}^{k}\sum_{j=0}^{l-1}\frac{\normMinv{\hatr_{l-1}}}{\normMinv{\hatr_j}}\biggr), \\
        \abs{\delta f_k} &\leq \bigO(nk\macheps)\hatalphak\norm{M^{-1}}^{\frac{1}{2}} \biggl(1+\sum_{j=0}^{k-1}\frac{\normMinv{\hatr_{k}}}{\normMinv{\hatr_j}}\biggr)\normMinv{\hatr_{k}}\norm{A}\max_{j\leq k+1}\bigl(\norm{\hatx_j},\norm{x}\bigr).
    \end{align}
\end{lemma}

\begin{proof}
By~\eqref{eq:tildexk1} and the definition of \(f\), \(f(\barxkone)-f(\barxk)\) can be written as
\begin{equation} \label{eq:fxk1-fx-0}
    \begin{split}
        f(\barxkone)-f(\barxk) &= f(\barxk+\hatalphak\hatpk+\Dx_{k+1})-f(\barxk) \\
        &= \frac{1}{2}(\barxk+\hatalphak\hatpk+\Dx_{k+1})\trans A(\barxk+\hatalphak\hatpk+\Dx_{k+1}) \\
        &\quad- b\trans (\barxk+\hatalphak\hatpk+\Dx_{k+1}) - \frac{1}{2} \normA{\barxk}^2 + b\trans \barxk \\
        &= \hatalphak\hatpk\trans A\barxk + \Dx_{k+1}\trans A\barxk + \frac{1}{2}\hatalphak^2 \normA{\hatpk}^2 + \hatalphak\hatpk\trans A\Dx_{k+1} \\
        &\quad+ \frac{1}{2} \normA{\Dx_{k+1}}^2 - b\trans \hatalphak\hatpk - b\trans\Dx_{k+1}.
    \end{split}
\end{equation}
Notice that, by~\eqref{eq:tildexk1},  \(\hatalphak\hatpk\trans A\barxkone\)  can be written as
\begin{align*}
\hatalphak\hatpk\trans A\barxkone &= \hatalphak\hatpk\trans A\barxk + \hatalphak^2 \normA{\hatpk}^2 + \hatalphak\hatpk\trans A\Dx_{k+1}.
\end{align*}
Together with the definition of \(\barrk\) and \(\barrkone\), \eqref{eq:fxk1-fx-0} can be simplified as
\begin{equation} \label{eq:fxk1-fxk-1}
    \begin{split}
        f(\barxkone)-f(\barxk) &= \hatalphak\hatpk\trans A\barxkone + \Dx_{k+1}\trans A\barxk - \frac{1}{2}\hatalphak^2 \normA{\hatpk}^2 \\
        &\quad+ \frac{1}{2} \normA{\Dx_{k+1}}^2 - b\trans \hatalphak\hatpk - b\trans\Dx_{k+1} \\
        &= -\hatalphak\hatpk\trans \barrkone - \Dx_{k+1}\trans \barrk
        - \frac{1}{2}\hatalphak^2\normA{\hatpk}^2 + \frac{1}{2} \normA{\Dx_{k+1}}^2.
    \end{split}
\end{equation}
Substituting the expression in~\eqref{eq:lem:normDalpha} for \(\hatalphak\) in~\eqref{eq:fxk1-fxk-1}, we have
\begin{equation} \label{eq:fxk1-fxk-2}
    \begin{split}
        f(\barxkone)-f(\barxk)
        &= -\hatalphak\hatpk\trans \barrkone - \Dx_{k+1}\trans \barrk
        - \frac{1}{2}\hatalphak(1+\Dalpha_k)\normMinv{\hatrk}^2 + \frac{1}{2} \normA{\Dx_{k+1}}^2.
    \end{split}
\end{equation}
Using~\eqref{eq:normDbarr}, we can substitute \(\barrk\) and \(\barrkone\) with \(\hatrk\) and \(\hatrkone\), respectively, in~\eqref{eq:fxk1-fxk-2}, i.e.,
\begin{equation} \label{eq:fxk1-fxk-3}
    \begin{split}
        f(\barxkone)-f(\barxk)
        &= - \frac{1}{2}\hatalphak(1+\Dalpha_k)\normMinv{\hatrk}^2 + \frac{1}{2}\normA{\Dx_{k+1}}^2
        -\hatalphak\hatpk\trans\hatrkone \\
        &\quad + \hatalphak\hatpk\trans\Delta\bar r_{k+1} - \Dx_{k+1}\trans \hatrk + \Dx_{k+1}\trans\Delta\bar r_{k} \\
        &= -\frac{1}{2}\hatalphak(1+\Dhatalpha_k)\normMinv{\hatrk}^2 + \delta f_k,
    \end{split}
\end{equation}
where \(\delta f_k := \hatalphak\hatpk\trans\Delta\bar r_{k+1}+ \Dx_{k+1}\trans\Delta\bar r_{k}\) and
\begin{equation*}
\begin{split}
    \Dhatalpha_k := \frac{\normA{\Dx_{k+1}}^2/2 -\hatalphak\hatpk\trans\hatrkone - \Dx_{k+1}\trans \hatrk}{-\hatalphak\normMinv{\hatrk}^2/2}
    = \frac{-\normA{\Dx_{k+1}}^2 + 2\hatalphak\hatpk\trans\hatrkone + 2\Dx_{k+1}\trans \hatrk}{\hatalphak\normMinv{\hatrk}^2}.
\end{split}
\end{equation*}
Then we will estimate the terms involved in \(\Dhatalpha_k\) and \(\delta f_k\) separately.
From Lemma~\ref{lem:rkpk-rk1pk}, \(\hatalphak\hatpk\trans\hatrkone\) satisfies
\begin{equation} \label{eq:lem:proof:fxk1-fxk-3-right-1}
    \begin{split}
        \abs{\hatalphak\hatpk\trans\hatrkone}
        &\leq \biggl(\bigO\bigl(nk\macheps\bigr)\kappa(A)+ \bigO\bigl(k\epspre\bigr)\biggr)\hatalphak\normMinv{\hatr_{k}}^2 \\
        &\quad+ \bigO\bigl(\macheps\bigr)\kappa(M)^{\frac{1}{2}} 
        \biggl(k+1+\sum_{l=0}^{k}\sum_{j=0}^{l-1}\frac{\normMinv{\hatr_{l-1}}}{\normMinv{\hatr_j}}\biggr)\hatalphak\normMinv{\hatr_{k}}^2.
    \end{split}
\end{equation}
From~\eqref{eq:normDbarr} and~\eqref{eq:xk}, we have
\begin{equation*} 
    \begin{split}
        \abs{\hatalphak\hatpk\trans\Delta\bar r_{k+1}}
        &\leq \hatalphak\norm{\hatpk}\norm{\Delta\bar r_{k+1}}
        \leq \bigO(nk\macheps)\hatalphak\norm{\hatpk}\norm{A}\max_{j\leq k+1}\bigl(\norm{\hatx_j},\norm{x}\bigr) \quad\text{and} \\
        \abs{\Dx_{k+1}\trans\Delta\bar r_{k}} &\leq \bigO(\macheps)\hatalphak\norm{\hatpk}\norm{\Delta\bar r_{k}}
        \leq \bigO(nk\macheps^2)\hatalphak\norm{\hatpk}\norm{A}\max_{j\leq k+1}\bigl(\norm{\hatx_j},\norm{x}\bigr).
    \end{split}
\end{equation*}
Then substituting~\eqref{eq:lem:u-normpk} for \(\macheps\norm{\hatpk}\) above, we have
\begin{equation} \label{eq:lem:proof:fxk1-fxk-3-right-2}
    \begin{split}
        \abs{\hatalphak\hatpk\trans\Delta\bar r_{k+1}}
        &\leq \bigO(nk\macheps)\hatalphak\norm{M^{-1}}^{\frac{1}{2}} \biggl(1+\sum_{j=0}^{k-1}\frac{\normMinv{\hatr_{k}}}{\normMinv{\hatr_j}}\biggr)\normMinv{\hatr_{k}}
        \norm{A}\max_{j\leq k+1}\bigl(\norm{\hatx_j},\norm{x}\bigr), \\
        \abs{\Dx_{k+1}\trans\Delta\bar r_{k}} &\leq \bigO(nk\macheps^2)\hatalphak\norm{M^{-1}}^{\frac{1}{2}} \biggl(1+\sum_{j=0}^{k-1}\frac{\normMinv{\hatr_{k}}}{\normMinv{\hatr_j}}\biggr)\normMinv{\hatr_{k}}
        \norm{A}\max_{j\leq k}\bigl(\norm{\hatx_j},\norm{x}\bigr).
    \end{split}
\end{equation}
By~\eqref{eq:xk} and~\eqref{eq:lem:u-normpk}, \(\Dx_{k+1}\trans \hatrk\) can be bounded as
\begin{equation} \label{eq:lem:proof:fxk1-fxk-3-right-3}
    \begin{split}
        \abs{\Dx_{k+1}\trans \hatrk}
        &\leq \bigO(\macheps)\hatalphak\norm{\hatpk}\norm{\hatrk} \\
        &\leq \bigO(\macheps)\hatalphak\norm{\hatpk}\bigl(\norm{M}^{\frac{1}{2}}\normMinv{\hatrk}\bigr) \\
        &\leq \bigO(\macheps)\kappa(M)^{\frac{1}{2}} 
        \biggl(1+\sum_{j=0}^{k-1}\frac{\normMinv{\hatrk}}{\normMinv{\hatr_j}}\biggr)\hatalphak\normMinv{\hatrk}^2.
    \end{split}
\end{equation}
From \eqref{eq:lem:proof:normAnormpk2/pkApk}, together with the assumptions \(\bigO(n\macheps)\kappa(A)\leq 1/2\) and \(\bigO\bigl(k\epspre\bigr) + \bigO(nk\macheps)\leq 1/2\), and observing that
\[
\abs{1+\Dalpha_k}
\leq 1+\abs{\Dalpha_k}
\leq 1 + \bigO(n\macheps)\kappa(A) + \bigO(\epspre)
\leq 2,
\]
we can bound \( \normA{\Dx_{k+1}}^2\) as
\begin{equation} \label{eq:lem:proof:fxk1-fxk-3-right-5}
    \normA{\Dx_{k+1}}^2
    \leq \bigO(\macheps^2) \hatalphak^2\norm{A}\norm{\hatpk}^2
    \leq \bigO(\macheps) \hatalphak\normMinv{\hatrk}^2.
\end{equation}
Combining~\eqref{eq:fxk1-fxk-3} with~\eqref{eq:lem:proof:fxk1-fxk-3-right-1}--\eqref{eq:lem:proof:fxk1-fxk-3-right-5}, we conclude the proof.
\end{proof}

In Lemma~\ref{lem:f-normrk}, we have established the relationship between \(f(\barxkone)-f(\barxk)\) and \(\normMinv{\hatrk}\), namely, that \(f(\barxkone)-f(\barxk)\) can be expressed as a quadratic function of \(\normMinv{\hatrk}\). 
Subsequently, we will demonstrate that the difference \(f(\barxk)-f(\barxkone)\) can be made sufficiently small by increasing \(k\), indicating that the value of the quadratic function of \(\normMinv{\hatrk}\) can be arbitrarily close to \(0\).

\begin{lemma} \label{lem:Dfsmall}
    Assume that \(\barxk\) and \(\barxkone\) satisfy~\eqref{eq:tildexk1}, and \(\hatalphak\) satisfies~\eqref{eq:alphak}.
    Let \(\epspre\) be as defined in Lemma~\ref{lem:eachstep-alphabeta}.
    If \(\bigO(n\macheps)\kappa(A)+\bigO(\epspre)\leq 1/2\), then for any \(\epsilon\geq 0\), there exists an iteration step \(\kstar\geq 0\) such that
    \[
        f(\barx_{\kstar}) - f(\barx_{\kstar+1}) \leq \epsilon\cdot c(\kstar)\cdot\hatalpha_{\kstar},
    \]
    and
    \begin{equation} \label{eq:lem-Dfsmall:sumrk/rj}
        \sum_{j=0}^{\kstar-1}\frac{\normMinv{\hatr_{\kstar}}}{\normMinv{\hatr_j}}\leq 3\kstar\kappa\bigl(M_L^{-1}AM_R^{-1}\bigr)^{1/2},
    \end{equation}
    where \(c(\kstar) := (\kstar)^2\biggl(1+\biggl(1+\bigO(\epspre)+ \bigO(n\macheps)\biggr)^{\kstar}\biggl(\bigO\bigl(\kstar\epspre\bigr) + \bigO\bigl(n\kstar\macheps\bigr)\biggr)\biggr)\).
\end{lemma}

\begin{proof}
This proof includes three steps:
\begin{enumerate}
    \item First, we will show that
    \begin{equation} \label{eq:lem-Dfsmall:proof:sumrk/rj}
        \sum_{j=0}^{k-1}\frac{\normMinv{\hatrk}}{\normMinv{\hatr_j}}\leq 3k\kappa\bigl(M_L^{-1}AM_R^{-1}\bigr)^{1/2}
    \end{equation}
    if \(f(\barxk)\leq f(\barx_j)\) holds for any \(j<k\).\\
    \item Subsequently, we will prove that for any \(\epsilon\geq 0\), there exists \(\kstar\geq 0\) such that \(f(\barx_{\kstar}) - f(\barx_{\kstar+1}) \leq \epsilon\).
    Moreover, for this specific \(\kstar\), \(f(\barx_{\kstar})\leq f(\barx_j)\) holds for any \(j<\kstar\).
    The existence of \(\kstar\) is in fact based on the common situation where \(f(\bar{x}_{\kstar}) < f(\bar{x}_{\kstar+1})\) once rounding error effects become dominant.\\
    \item Finally, we will conclude by deriving a lower bound for \(\hatalpha_{\kstar}\).\\
\end{enumerate}

We begin with proving~\eqref{eq:lem-Dfsmall:proof:sumrk/rj}, which requires a bound on \(\normMinv{\hatrk}/\normMinv{\hatr_j}\) for any \(j<k\).
Using~\eqref{eq:hatrk-barrk} we can express \(\normMinv{\hatrk}\) as a perturbation of $\normMinv{\barrk}$, i.e.,
\begin{equation}
    \begin{split}
        \normMinv{\hatrk} &= \normMinv{\barrk}(1 + \delta g_{k}),
    \end{split}
\end{equation}
with the relative error \(\delta g_{k}\) satisfying
\begin{equation} \label{eq:lem-proof:deltagk}
    \begin{split}
        \abs{\delta g_{k}} &= \frac{\abs{\normMinv{\hatrk}-\normMinv{\barrk}}}{\normMinv{\barrk}} \\
        &\leq \frac{\normMinv{\overbrace{\hatrk-\barrk}^{= \Delta\barrk}}}{\normMinv{\barrk}} \\
        &\leq \frac{\bigO(nk\macheps) \norm{A} \max_{j\leq k}\bigl(\norm{\hatx_j}, \norm{x}\bigr)\norm{M^{-1}}^{\frac{1}{2}}}{\normMinv{\barrk}},
    \end{split}
\end{equation}
where the last inequality follows from~\eqref{eq:normDbarr}.
Notice that either \(\abs{\delta g_{k}}\leq 1/2\) or \(\abs{\delta g_{k}}> 1/2\) holds.
Here \(1/2\) can be replaced by any arbitrary number \(1/d\) with \(d=\bigO(1)\).
In the second case, we can estimate the denominator \(\normMinv{\barrk}\) in~\eqref{eq:lem-proof:deltagk} by \(\norm{\barrk}\leq \norm{M}^{1/2}\normMinv{\barrk}\), which gives that
\begin{equation} \label{eq:lem-proof:normbarrk-normMinvbarrk}
    \norm{\barrk}\leq \norm{M}^{1/2}\normMinv{\barrk}
    \leq \bigO(nk\macheps) \kappa(M^{-1})^{\frac{1}{2}} \norm{A} \max_{j\leq k}\bigl(\norm{\hatx_j}, \norm{x}\bigr).
\end{equation}
Then, using the results from Theorem~\ref{thm:Drk}, together with~\eqref{eq:hatxk1-tildexk1} and~\eqref{eq:lem-proof:normbarrk-normMinvbarrk}, we can bound the norm of the true residual as
\begin{equation}
    \norm{\tilde{r}_k}
    \leq \norm{\barrk} + \macheps k \norm{A} \max_{j \leq k} \norm{\hat{x}_j}
    \leq \bigO(nk\macheps) \kappa(M^{-1})^{\frac{1}{2}} \norm{A} \max_{j\leq k}\bigl(\norm{\hatx_j}, \norm{x}\bigr).
\end{equation}
The last equation gives us the backward error result for the PCG algorithm that we want to show.
Thus, without loss of generality, we can assume that \(\abs{\delta g_{k}}\leq 1/2\) here.
The previous analysis also holds for \(\normMinv{\hatr_j}\) and \(\delta g_j\) with \(j<k\).
Therefore, it follows that
\begin{equation*}
    \begin{split}
        \frac{\normMinv{\hatrk}}{\normMinv{\hatr_j}}
        &= \frac{\normMinv{\barrk}(1+\delta g_{k})}{\normMinv{\barr_{j}}(1+\delta g_{j})} 
        \leq \frac{3\normMinv{\barrk}}{\normMinv{\barr_{j}}} 
        = 3 \frac{\normMinv{A (x - \barxk)}}{\normMinv{A (x - \bar{x}_j)}}.
    \end{split}
\end{equation*}
Now, suppose that \(f(\barxk)\leq f(\barx_j)\) holds for any \(j<k\).
Using \eqref{eq:energy-min} we obtain
\[
    \frac{1}{2}\normA{\barxk-x}^2=f(\barxk)-f(x)\leq f(\barx_j)-f(x)=\frac{1}{2}\normA{\barx_j-x}^2,
\]
which implies that \(\normA{\barxk-x}\leq \normA{\barx_j-x}\) holds for any \(j<k\).
Furthermore, notice that for any $1 \leq i \leq k$ we have
\begin{align*}
    \normMinv{A(x - \bar{x}_i)} &= [(x - \bar{x}_i)^\top A M^{-1} A(x - \bar{x}_i)]^{\frac{1}{2}}\\
    &= [(x - \bar{x}_i)^\top A^\frac{1}{2} A^\frac{1}{2} M^{-1} A^\frac{1}{2} A^\frac{1}{2} (x - \bar{x}_i)]^{\frac{1}{2}}\\
    &\leq \norm{A^\frac{1}{2} M^{-1} A^\frac{1}{2}}^\frac{1}{2} \normA{x - \bar{x}_i}.
\end{align*}
Equivalently, it holds that 
\begin{align*}
    {\norm{A(x - \bar{x}_i)}}_{M^{-1 }}^{-1} 
    &\geq \frac{\norm{A^{-\frac{1}{2}} M A^{-\frac{1}{2}}} }{\normA{x - \bar{x}_i}}^\frac{1}{2}.
\end{align*}
Hence, in total, we have
\begin{equation} \label{eq:lem:proof:mnormrk1/mnormrkj}
    \frac{\normMinv{\hatrk}}{\normMinv{\hatr_j}}
    \leq 3\kappa\bigl(A^{\frac{1}{2}}M^{-1}A^{\frac{1}{2}}\bigr)^{\frac{1}{2}}\frac{\normA{\barxk-x}}{\normA{\barx_j-x}}
    \leq 3\kappa\bigl(A^{\frac{1}{2}}M^{-1}A^{\frac{1}{2}}\bigr)^{\frac{1}{2}}.
\end{equation}
Note that \(A\) and \(M^{-1}=M_R^{-1}M_L^{-1}\) are both symmetric positive definite matrices. Furthermore, since the matrices \(A^{\frac{1}{2}}M^{-1}A^{\frac{1}{2}}\) and \(M_L^{-1}AM_R^{-1} = \bigl(A^{\frac{1}{2}}M_R^{-1}\bigr)^{-1}A^{\frac{1}{2}}M^{-1}A^{\frac{1}{2}}\bigl(A^{\frac{1}{2}}M_R^{-1}\bigr)\) are similar, we further have
\begin{equation}
    \kappa\bigl(A^{\frac{1}{2}}M^{-1}A^{\frac{1}{2}}\bigr)
    = \frac{\lamdmax\bigl(M_L^{-1}AM_R^{-1}\bigr)}{\lamdmin\bigl(M_L^{-1}AM_R^{-1}\bigr)}
    = \kappa\bigl(M_L^{-1}AM_R^{-1}\bigr),
\end{equation}
which proves~\eqref{eq:lem-Dfsmall:proof:sumrk/rj} together with~\eqref{eq:lem:proof:mnormrk1/mnormrkj}.

We now move to the second step to prove that for any \(\epsilon\geq 0\), there exists \(\kstar\geq 0\) such that \(f(\barx_{\kstar}) - f(\barx_{\kstar+1}) \leq \epsilon\).
Here we consider two cases.\\
\begin{flushleft}
\underline{Case 1.} There exists \(\kstar\geq 0\) such that
\begin{equation*}
    f(\barx_{\kstar+1})-f(\barx_{\kstar}) \geq 0,
\end{equation*}
which implies that \(f(\barx_{\kstar}) - f(\barx_{\kstar+1})\leq 0\leq \epsilon\).\\
\end{flushleft}
\underline{Case 2.} The sequence \(\{f_k\}_{k=0}^{\infty}\) with \(f_k = f(\barxk)\) is monotonically decreasing.
Since \(f\) attains its minimum for the solution of the linear algebraic system \eqref{problem:linear-sys}, we have 
\[
f(x)\leq f(\barx_k)\leq f(x_0)
\]
for any \(k\geq 0\). Hence, \(\{f_k\}_{k=0}^{\infty}\) is a bounded sequence, and it is therefore convergent. That is, there exists \(f_\star\) such that \(\lim_{k\rightarrow +\infty} f(\barx_k) = \lim_{k\rightarrow +\infty} f_k = f_\star\).
Furthermore, this implies the existence of an iteration \(\kstar\geq 0\) such that \(f(\barx_{\kstar}) - f(\barx_{\kstar+1}) \leq \epsilon\).\\
Note that, for both these two cases, \(f(\barx_{\kstar})\leq f(\barx_j)\) holds for any \(j<\kstar\).
Taking the sum of~\eqref{eq:lem:proof:mnormrk1/mnormrkj} over \(j=0,\dotsc,\kstar\) yields~\eqref{eq:lem-Dfsmall:sumrk/rj}.

Now, we prove the final result.
Notice that for any \(k\), using the relations~\eqref{eq:lem:normDalpha}, \eqref{eq:lem:normpk-0}, and~\eqref{eq:lem-Dfsmall:sumrk/rj}, it holds that
\begin{equation} \label{eq:thm:proof:alphalowerbound}
\begin{split}
    \hatalphak
    &= \bigl(1+\Dalpha_k\bigr)\frac{\normMinv{\hatrk}^2}{\norm{\hatpk}^2} \frac{\norm{\hatpk}^2}{\hatpk\trans A\hatpk} \\
    &\geq \biggl(1-\bigO(n\macheps)\kappa(A) - \bigO(\epspre)\biggr)
    \frac{\normMinv{\hatrk}^2}{\norm{\hatpk}^2\norm{A}} \\
    &\geq \frac{1-\bigO(n\macheps)\kappa(A) - \bigO(\epspre)}
    {9\cdot c(k)\bigl(1+\kappa\bigl(M_L^{-1}AM_R^{-1}\bigr)^{\frac{1}{2}}\bigr)^2\norm{A}\norm{M^{-1}}},
\end{split}
\end{equation}
where \(c(k)\) has been defined in the statement of this lemma.
In the second step, we have proved that for any \(\epsilon \geq 0\), there exists \(\kstar \geq 0\) such that \(f(\barx_{\kstar}) - f(\barx_{\kstar+1}) \leq \epsilon\), where \(\epsilon\) cannot be chosen as a function of \(\kstar\).
From this conclusion, there exists \(\kstar\geq 0\) such that
\begin{equation}
\begin{split}
    f(\barx_{\kstar}) - f(\barx_{\kstar+1}) \leq 
    \frac{\epsilon \bigl(1-\bigO(n\macheps)\kappa(A)-\bigO(\epspre)\bigr)}{9\bigl(1+\kappa\bigl(M_L^{-1}AM_R^{-1}\bigr)^{\frac{1}{2}}\bigr)^2\norm{A}\norm{M^{-1}}}
    \leq \epsilon\cdot c(\kstar)\cdot\hatalpha_{\kstar},
\end{split}
\end{equation}
where the last inequality is derived from~\eqref{eq:thm:proof:alphalowerbound}.
\end{proof}

Lemma~\ref{lem:Dfsmall} proves the existence of the iteration \(\kstar\) and further provides the upper bound for \(\sum_{j=0}^{\kstar-1}\bigl(\normMinv{\hatr_{\kstar}}/\normMinv{\hatr_j}\bigr)\) at this specific iteration.
Combining the results of Lemma~\ref{lem:f-normrk} and Lemma~\ref{lem:Dfsmall}, we obtain a quadratic inequality with respect to \(\normMinv{\hatrk}\).
Solving this inequality provides an initial estimate for \(\normMinv{\hatrk}\).

\begin{lemma} \label{lem:backward-forward-err}
    Assume that \(\hatxk\) and \(\hatrk\) satisfy~\eqref{eq:xk} and~\eqref{eq:hatrk1}, and \(\barxk\) satisfies~\eqref{eq:tildexk1}.
    Let \(\epspre\) be as defined in Lemma~\ref{lem:eachstep-alphabeta}.
    If \(\bigO(n\macheps)\kappa(A)+\bigO\bigl(\epspre\bigr)\leq 1/2\), then there exists an iteration \(\kstar\geq 0\) such that
    \[
        f(\barx_{\kstar}) - f(\barx_{\kstar+1}) \leq c(\kstar)\cdot\bigO\bigl(n^2\macheps^2\bigr)\hatalpha_{\kstar}\norm{A}^2\norm{M^{-1}}\norm{x}^2,
    \]
    where \(c(\kstar)\) has been defined in Lemma~\ref{lem:Dfsmall}.
    Furthermore, if \(\kstar\) satisfies
    \begin{equation}  \label{eq:lem:backward-forward-err:assump}
    \begin{split}
        &\bigO\bigl(n\kstar\macheps\bigr)\kappa(A) + \bigO(\kstar\epspreq)
        + \bigO\bigl((\kstar)^2\macheps) \kappa(M)^{\frac{1}{2}} \kappa\bigl(M_L^{-1}AM_R^{-1}\bigr)^{\frac{1}{2}}\leq \frac{1}{2},
    \end{split}
    \end{equation}
    then
    \begin{align}
        \frac{\normMinv{\hatr_{\kstar}}}{\norm{A}\norm{M^{-1}}^{\frac{1}{2}} \norm{x}} 
        &\leq \bigO\bigl(n\kstar\macheps\bigr) \max_{j\leq \kstar+1}\biggl(\frac{\norm{\hatx_j}}{\norm{x}}, 1\biggr) \biggl(1+\sum_{j=0}^{\kstar-1}\frac{\normMinv{\hatr_{\kstar}}}
        {\normMinv{\hatr_j}}\biggr) \label{eq:thm:normMinvhatr}
    \end{align}
    and
    \begin{equation}
        \frac{\norm{\hatr_{\kstar}}}{\norm{A}\norm{x}} 
        \leq \bigO\bigl(n\kstar\macheps\bigr) \kappa(M)^{\frac{1}{2}} \max_{j\leq \kstar+1}\biggl(\frac{\norm{\hatx_j}}{\norm{x}}, 1\biggr) \biggl(1+\sum_{j=0}^{\kstar-1}\frac{\normMinv{\hatr_{\kstar}}}
        {\normMinv{\hatr_j}}\biggr). \label{eq:thm:normhatr}
    \end{equation}
\end{lemma}

\begin{proof}
    From Lemma~\ref{lem:f-normrk} we have 
    \[ 
        0 = - \frac{1}{2} \hatalphak (1 + \delta \hatalphak) \normMinv{\hatrk}^2 + \delta f_k + f(\barxk) - f(\barxkone).
    \]
    In particular, if we set \(\epsilon = \bigO(n^2\macheps^2)\norm{A}^2\norm{M^{-1}} \norm{x}^2\) in Lemma~\ref{lem:Dfsmall}, then there exists an iteration \(\kstar \geq 0 \) such that
    \begin{equation*}
        \begin{split}
            0&\leq - \frac{1}{2}\hatalpha_{\kstar}\biggl(1-\bigO\bigl(n\kstar\macheps\bigr)\kappa(A) - \bigO(\kstar\epspreq)
            - \bigO\bigl((\kstar)^2\macheps) \kappa(M)^{\frac{1}{2}} \kappa\bigl(M_L^{-1}AM_R^{-1}\bigr)^{\frac{1}{2}}\biggr) \normMinv{\hatr_{\kstar}}^2 \\
            &\quad +\bigO\bigl(n\kstar\macheps\bigr)\hatalpha_{\kstar} \norm{M^{-1}}^{\frac{1}{2}} \biggl(1+\sum_{j=0}^{\kstar-1}\frac{\normMinv{\hatr_{\kstar}}}{\normMinv{\hatr_j}}\biggr)\normMinv{\hatr_{\kstar}}\norm{A}\max_{j\leq \kstar+1}\bigl(\norm{\hatx_j},\norm{x}\bigr) \\
            &\quad + c(\kstar)\cdot\bigO\bigl(n^2\macheps^2\bigr) \hatalpha_{\kstar}\norm{A}^2\norm{M^{-1}} \norm{x}^2.
        \end{split}
    \end{equation*}
    To simplify the notation for the resulting quadratic inequality above, we introduce the coefficients \(\gamma_a\), \(\gamma_b\), and \(\gamma_c\) such that
    \begin{equation*}
        -\gamma_a \normMinv{\hatr_{\kstar}}^2 + \gamma_b\normMinv{\hatr_{\kstar}} + \gamma_c\geq 0
    \end{equation*}
    with
    \begin{align*}
        \gamma_a &= \frac{1}{2}\hatalpha_{\kstar}\biggl(1-\bigO\bigl(n\kstar\macheps\bigr)\kappa(A) - \bigO(\kstar\epspreq)
        - \bigO\bigl((\kstar)^2\macheps) \kappa(M)^{\frac{1}{2}} \kappa\bigl(M_L^{-1}AM_R^{-1}\bigr)^{\frac{1}{2}}\biggr), \\
        \gamma_b &= \bigO\bigl(n\kstar\macheps\bigr)\hatalpha_{\kstar} \norm{M^{-1}}^{\frac{1}{2}} \biggl(1+\sum_{j=0}^{\kstar-1}\frac{\normMinv{\hatr_{\kstar}}}{\normMinv{\hatr_j}}\biggr)\norm{A}\max_{j\leq \kstar+1}\bigl(\norm{\hatx_j},\norm{x}\bigr), \\
        \gamma_c &= c(\kstar)\cdot\bigO\bigl(n^2\macheps^2\bigr) \hatalpha_{\kstar}\norm{A}^2\norm{M^{-1}} \norm{x}^2.
    \end{align*}
    Notice that \(\gamma_a>0\) due to the assumption~\eqref{eq:lem:backward-forward-err:assump}, which means that \(\gamma_a\) can safely be used as a denominator.
    By solving this quadratic inequality, we then have
    \begin{equation} \label{eq:thm:proof:rk-beta}
    \begin{split}
        \normMinv{\hatr_{\kstar}} \leq \frac{-\gamma_b - \sqrt{\gamma_b^2+4\gamma_a\gamma_c}}{-2\gamma_a}
        \leq \frac{\gamma_b}{\gamma_a} + \frac{\sqrt{\gamma_c}}{\sqrt{\gamma_a}}.
    \end{split}
    \end{equation}
    Due to assumption~\eqref{eq:lem:backward-forward-err:assump}, the term \(\gamma_b/\gamma_a\) in the right-hand side of~\eqref{eq:thm:proof:rk-beta} satisfies
    \begin{equation} \label{eq:thm:proof:rkbeta-bound1}
        \begin{split}
            \frac{\gamma_b}{\gamma_a} &\leq \bigO\bigl(n\kstar\macheps\bigr)\norm{A}\max_{j\leq \kstar+1}\bigl(\norm{\hatx_j}, \norm{x}\bigr) \norm{M^{-1}}^{\frac{1}{2}} 
            \biggl(1+\sum_{j=0}^{\kstar-1}\frac{\normMinv{\hatr_{\kstar}}}{\normMinv{\hatr_j}}\biggr).
        \end{split}
    \end{equation}
    Similarly, the term \(\sqrt{\gamma_c}\,/\sqrt{\gamma_a}\) satisfies
    \begin{equation} \label{eq:thm:proof:rkbeta-bound2}
    \begin{split}
        \frac{\sqrt{\gamma_c}}{\sqrt{\gamma_a}} &\leq \sqrt{\frac{c(\kstar)\cdot\bigO\bigl(n^2\macheps^2\bigr) \hatalpha_{\kstar}\norm{A}^2\norm{M^{-1}} \norm{x}^2}
        {\hatalpha_{\kstar}/4}} \\
        &\leq \sqrt{c(\kstar)}\cdot\bigO(n\macheps)\norm{A}\norm{M^{-1}}^{\frac{1}{2}} \norm{x} \\
        &\leq \bigO(n\kstar\macheps)\norm{A}\norm{M^{-1}}^{\frac{1}{2}} \norm{x}.
    \end{split}
    \end{equation}
    Substituting \(\gamma_b/\gamma_a\) and \(\sqrt{\gamma_c}\,/\sqrt{\gamma_a}\) in~\eqref{eq:thm:proof:rk-beta} by~\eqref{eq:thm:proof:rkbeta-bound1} and~\eqref{eq:thm:proof:rkbeta-bound2}, respectively, we derive~\eqref{eq:thm:normMinvhatr}, i.e., the bound on \(\normMinv{\hatr_{\kstar}}\).

    Furthermore, since
    \[
        \lambda_{\min}(M^{1/2}) \normMinv{\hatr_{\kstar}} \leq \norm{\hatr_{\kstar}} \leq \norm{M^{1/2}} \norm{M^{-1/2}\hatr_{\kstar}},
    \]
    we obtain~\eqref{eq:thm:normhatr} due to \(\kappa(M)^{1/2}=\norm{M^{1/2}}/\lamdmin(M^{1/2})=\kappa(M^{-1})^{1/2}\).
\end{proof}

Within the bound on \(\normMinv{\hatr_{\kstar}}\) shown in Lemma~\ref{lem:backward-forward-err}, the term \(\sum_{j=0}^{\kstar-1}\bigl(\normMinv{\hatr_{\kstar}}/\normMinv{\hatr_j}\bigr)\) is the only unknown factor.
It is vital to emphasize that~\eqref{eq:lem-Dfsmall:sumrk/rj} provides an upper bound for \(\sum_{j=0}^{\kstar-1}\bigl(\normMinv{\hatr_{\kstar}}/\normMinv{\hatr_j}\bigr)\) associated with \(\kappa\bigl(M_L^{-1}AM_R^{-1}\bigr)^{\frac{1}{2}}\).
However, in practice, \(\norm{\hatrk}\) can always converge to a relatively small level as long as \(A\) and \(M\) are not too ill-conditioned, suggesting that the bound for \(\normMinv{\hatrk}\) should not involve \(\kappa\bigl(M_L^{-1}AM_R^{-1}\bigr)^{\frac{1}{2}}\), especially for the case \(M_L=M_R=I\).
Therefore, our aim here is to demonstrate that, for a certain iteration, \(\sum_{j=0}^{\kstar-1}\bigl(\normMinv{\hatr_{\kstar}}/\normMinv{\hatr_j}\bigr)\) can be \(\bigO(\kstar)\) without being tied to \(\kappa\bigl(M_L^{-1}AM_R^{-1}\bigr)^{\frac{1}{2}}\).
Furthermore, according to Lemma~\ref{lem:backward-forward-err}, \(\norm{\hatrk}\) can be bounded.

\begin{theorem} \label{thm:backward-forward-err}
    Assume that \(\hatxk\) and \(\hatrk\) are generated by Algorithm~\ref{alg:pcg}, and \(\barxk\) satisfies~\eqref{eq:tildexk1}.
    Let \(\epspre\) be as defined in Lemma~\ref{lem:eachstep-alphabeta}.
    If \(\bigO(n\macheps)\kappa(A)+\bigO\bigl(\epspre\bigr)\leq 1/2\), then there exists an iteration step \(\kstar\geq 0\) such that
    \[
        f(\barx_{\kstar}) - f(\barx_{\kstar+1}) \leq c(\kstar)\cdot\bigO\bigl(n^2\macheps^2\bigr)\hatalpha_{\kstar}\norm{A}^2\norm{M^{-1}} \norm{x}^2,
    \]
    where \(c(\kstar)\) has been defined in Lemma~\ref{lem:Dfsmall}.
    Furthermore, if \(\kstar\) satisfies~\eqref{eq:lem:backward-forward-err:assump},
    then there exists \(i\leq \kstar\) such that
    \begin{align}
        \frac{\norm{\hatr_{i}}}{\norm{A}\norm{x}}&\leq \bigO\bigl(n(\kstar)^2\macheps\bigr) \kappa(M)^{\frac{1}{2}}
        \max_{j\leq \kstar+1}\biggl(\frac{\norm{\hatx_j}}{\norm{x}}, 1\biggr), \label{eq:thm-final:normhatr} \\
        \frac{\norm{b-A\hatx_{i}}}{\norm{A}\norm{x}}&\leq \bigO\bigl(n(\kstar)^2\macheps\bigr) \kappa(M)^{\frac{1}{2}}
        \max_{j\leq \kstar+1}\biggl(\frac{\norm{\hatx_j}}{\norm{x}}, 1\biggr), \label{eq:thm-final:normtilder}
    \end{align}
    and
    \begin{equation} \label{eq:thm-final:normA-xk-x}
        \begin{split}
            \frac{\normA{\hatx_{i}-x}}{\norm{A}^{1/2}\norm{x}}
            &\leq \bigO\bigl(n(\kstar)^2\macheps\bigr)\kappa(M)^{\frac{1}{2}} \kappa(A)^{1/2} \max_{j\leq \kstar+1}\biggl(\frac{\norm{\hatx_j}}{\norm{x}}, 1\biggr).
        \end{split}
    \end{equation}
\end{theorem}

\begin{proof}
    From Lemma~\ref{lem:backward-forward-err}, there exists \(\kstar\geq 0\) such that
    \begin{equation}
        f(\barx_{\kstar}) - f(\barx_{\kstar+1}) \leq c(\kstar)\cdot\bigO\bigl(n^2\macheps^2\bigr)\hatalpha_{\kstar}\norm{A}^2 \norm{M^{-1}}\norm{x}^2
    \end{equation}
    and~\eqref{eq:thm:normhatr} hold.

    First, we aim to prove~\eqref{eq:thm-final:normhatr}.
    From Lemma~\ref{lem:backward-forward-err}, we know that the upper bound on \(\normMinv{\hatr_{\kstar}}\) is influenced mainly by \(\sum_{j=0}^{\kstar-1}\bigl(\normMinv{\hatr_{\kstar}}/\normMinv{\hatr_j}\bigr)\).
    Thus, we only need to estimate the term \(\sum_{j=0}^{\kstar-1}\bigl(\normMinv{\hatr_{\kstar}}/\normMinv{\hatr_j}\bigr)\).
    For the estimate of \(\normMinv{\hatr_{\kstar}}/\normMinv{\hatr_j}\), we make the following case distinction:\\
    \begin{flushleft}     
    \underline{Case 1.} It holds that \(\normMinv{\hatr_{\kstar}}\leq \normMinv{\hatr_j}\) for all \(j\leq \kstar\).
    For this case, we have
    \begin{equation} \label{eq:thm:proof:normpk/normrk}
        \sum_{j=0}^{\kstar-1}\frac{\normMinv{\hatr_{\kstar}}}{\normMinv{\hatr_j}} \leq \kstar.
    \end{equation}
    Together with~\eqref{eq:thm:normhatr} in Lemma~\ref{lem:backward-forward-err}, we obtain~\eqref{eq:thm-final:normhatr} with \(i=\kstar\), i.e., the bound on \(\norm{\hatr_i}\).\\
    \end{flushleft}
\begin{flushleft}
    
    \underline{Case 2.} It does not hold that \(\normMinv{\hatr_{\kstar}}\leq \normMinv{\hatr_j}\) for all \(j\leq \kstar\).
    Thus, there exists an iteration \(i < \kstar \) such that
    \[
    \min_{1 \leq k \leq \kstar} (\normMinv{\hatrk}) = \normMinv{\hat{r}_i} = \tilde c \normMinv{\hatr_{\kstar}}
    \]
    with \(\tilde c < 1\).
    This implies that~\eqref{eq:thm:proof:normpk/normrk} also holds for \(i\), i.e.,
    \begin{equation} \label{eq:thm:proof:normpi/normri}
        \sum_{j=0}^{i-1}\frac{\normMinv{\hatr_{i}}}{\normMinv{\hatr_j}}\leq i.
    \end{equation}
    In particular,
    \begin{equation*}
    \begin{split}
        \sum_{j = 0}^{\kstar - 1} \frac{\normMinv{\hatr_\kstar}}{\normMinv{\hatr_j}}
        =\frac{1}{\tilde c} \biggl(\sum_{j = 0}^{i - 1} \frac{\normMinv{\hatr_i}}{\normMinv{\hatr_j}} + \sum_{l = i}^{\kstar - 1} \frac{\normMinv{\hatr_i}}{\normMinv{\hatr_l}} \biggr)
        \leq \frac{\kstar}{\tilde c}.
    \end{split}
    \end{equation*}
    Together with Lemma~\ref{lem:backward-forward-err}, \(\normMinv{\hatr_{\kstar}}\) satisfies
    \begin{equation} \label{eq:thm:proof:rkstarbound}
    \begin{split}
        \normMinv{\hatr_{\kstar}}
        &\leq \bigO\bigl(n\kstar\macheps\bigr) \norm{A} \norm{M^{-1}}^{\frac{1}{2}}
        \max_{j\leq \kstar+1}\bigl(\norm{\hatx_j}, \norm{x}\bigr)
        \biggl(1+\sum_{j=0}^{\kstar-1}\frac{\normMinv{\hatr_{\kstar}}}{\normMinv{\hatr_j}}\biggr)\\
         &\leq \bigO\bigl(n(\kstar)^2\macheps\bigr) 
        \frac{\norm{A}}{\tilde c} \norm{M^{-1}}^{\frac{1}{2}}
        \max_{j\leq \kstar+1}\bigl(\norm{\hatx_j}, \norm{x}\bigr).
    \end{split}
    \end{equation}
\end{flushleft}
    We now aim to prove \(f(\barx_i) - f(\barx_{i+1})\) is sufficiently small for this case.
    From the assumption~\eqref{eq:lem:backward-forward-err:assump}, we have \(\abs{\delta\hatalpha_i}\leq 1\) and thus \(\hatalpha_i(1+\Dhatalpha_i)\normMinv{\hat{r}_i}^2/2\leq \hatalpha_i\normMinv{\hatr_i}^2\).
    Together with Lemma~\ref{lem:f-normrk}, \eqref{eq:thm:proof:normpi/normri}, and \(\min_{1\leq k\leq\kstar}\bigl(\normMinv{\hatrk}\bigr) = \normMinv{\hatr_i} = \tilde c\normMinv{\hatr_{\kstar}}\) with \(0\leq \tilde c<1\), we have
    \begin{equation}
        \begin{split}
            f(\barx_i) - f(\barx_{i+1})
            &\leq \hatalpha_i\normMinv{\hatr_i}^2 
            + \bigO(ni^2\macheps)\hatalpha_i\norm{M^{-1}}^{\frac{1}{2}} \normMinv{\hatr_{i}}\norm{A}\max_{j\leq i+1}\bigl(\norm{\hatx_j},\norm{x}\bigr) \\
            &\leq \hatalpha_i \tilde c^2\normMinv{\hatr_{\kstar}}^2 
            + \bigO(ni^2\macheps)\hatalpha_i \tilde c\norm{M^{-1}}^{\frac{1}{2}}
            \normMinv{\hatr_{\kstar}}\norm{A}\max_{j\leq i+1}\bigl(\norm{\hatx_j},\norm{x}\bigr),
        \end{split}
    \end{equation}
    where the first inequality is obtained from the application of assumption~\eqref{eq:lem:backward-forward-err:assump} on the bound for \(\Dhatalpha_i\) in Lemma~\ref{lem:f-normrk}.
    Using the bound on \(\normMinv{\hatr_{\kstar}}\) shown in~\eqref{eq:thm:proof:rkstarbound}, we derive
    \begin{equation}
        \begin{split}
            f(\barx_i) - f(\barx_{i+1})
            &\leq \bigO(n^2(\kstar)^4\macheps^2) \hatalpha_i \norm{A}^2\norm{M^{-1}}\max_{j\leq \kstar+1}\bigl(\norm{\hatx_j},\norm{x}\bigr)^2,
        \end{split}
    \end{equation}
    which gives~\eqref{eq:thm-final:normhatr} by following the analysis of Lemma~\ref{lem:backward-forward-err}.

    Now it remains to prove~\eqref{eq:thm-final:normtilder} and~\eqref{eq:thm-final:normA-xk-x}.
    To do this, we use the relation of the true residual to the recursively updated residual in Theorem~\ref{thm:Drk} and the \(A\)-norm of the error, respectively.
    Combining~\eqref{eq:thm-final:normhatr} with Theorem~\ref{thm:Drk}, we further obtain~\eqref{eq:thm-final:normtilder}, i.e., the bound on \(\norm{\tilde r_{i}}\) with \(\tilde r_{i} = b-A\hatx_{i}\).
    Notice that, for any \(i\),
    \begin{equation} \label{eq:thm:proof:tilder-xk-x-A}
        \norm{\tilde r_i} = \norm{b-A\hatx_i} = \norm{A(x-\hatx_i)}
        \geq \lamdmin(A^{1/2})\norm{A^{1/2}(x-\hatx_i)}
        = \lamdmin(A^{1/2})\normA{\hatx_i-x}.
    \end{equation}
    Thus, \(\normA{\hatx_{i}-x}\) can be bounded as
    \begin{equation*}
        \begin{split}
            \normA{\hatx_{i}-x}&\leq \frac{\norm{\tilde r_{i}}}{\lamdmin(A)^{1/2}},
        \end{split}
    \end{equation*}
    which gives~\eqref{eq:thm-final:normA-xk-x} together with~\eqref{eq:thm-final:normtilder}.
\end{proof}

In contrast to other works such as~\cite{Greenbaum1997,ST2002,ST2005}, Theorem~\ref{thm:backward-forward-err} specifically addresses the backward error, rather than a bound on the residual gap under the assumption that the norm of the recursively updated residual goes to zero. Moreover, \eqref{eq:thm-final:normA-xk-x} also provides an upper bound for the forward error in the \(A\)-norm for Algorithm~\ref{alg:pcg}.

\begin{remark} \label{remark:kstar}
    In Theorem~\ref{thm:backward-forward-err}, there is a connection between the backward and forward errors and the number of iterations, suggesting that PCG may fail to reach \(\bigO(\macheps)\) level for the backward and forward errors if the convergence rate is excessively slow.
    However, to the best of our knowledge, it is hard to construct examples at which PCG fails when \(\norm{x_0}/\norm{x}=\bigO(1)\).
    This is largely due to the fact that the impact of the number of iterations on rounding errors is similar to that of the dimension, which tends to be considerably less significant in practical scenarios than theoretical predictions indicate.
    For example, in~\cite{CHM2021,ESDP2023,IZ2020}, the authors show that the stochastic rounding error bounds for the basic linear algebra kernels can be proportional to \(\sqrt{n}\macheps\) instead of \(n\macheps\).
    Moreover, PCG is typically used to solve large linear systems, hence the maximum number of iterations acceptable by users is generally much smaller than \(n\). 
    
    Furthermore, it is worth emphasizing that, although our analysis shows that the norm of the recursively updated residual can become very small---particularly in the unpreconditioned case---the bound on the right-hand side of~\eqref{eq:thm-final:normhatr} is clearly a significant overestimate.
    Indeed, as mentioned in the introduction, it is commonly observed in practice that the norm of the recursively updated residual continues to decrease long after the true residual has reached its maximal attainable accuracy.   
\end{remark}

\begin{remark}
    From~\eqref{eq:thm-final:normtilder} or~\eqref{eq:thm-final:normA-xk-x}, we can also derive the forward error in the \(2\)-norm, that is,
    \begin{equation}
        \frac{\norm{\hatx_{i}-x}}{\norm{x}}
        \leq \bigO\bigl(n(\kstar)^2\macheps\bigr)\kappa(A) \kappa(M^{-1})^{\frac{1}{2}} \max_{j\leq \kstar+1}\biggl(\frac{\norm{\hatx_j}}{\norm{x}}, 1\biggr),
    \end{equation}
    and obtain the standard relative backward error, that is,
    \begin{equation}
        \frac{\norm{b-A\hatx_{i}}}{\norm{A}\norm{\hatx_i}+\norm{b}} 
        \leq\frac{\norm{b-A\hatx_{i}}}{\norm{A}\norm{\hatx_i}}
        \leq\bigO\bigl(n(\kstar)^2\macheps\bigr) \kappa(M^{-1})^{\frac{1}{2}} 
        \max_{j\leq \kstar+1}\biggl(\frac{\norm{\hatx_j}}{\norm{\hatx_i}}, \frac{\norm{x}}{\norm{\hatx_i}}\biggr).
    \end{equation}
\end{remark}

\begin{remark}
    If we do not use any preconditioner in the CG algorithm, i.e., \(M_R=M_L=I\) and \(\epss=\epsq=\epsz=0\), then, from Theorem~\ref{thm:backward-forward-err}, the backward and forward errors can be bounded as
    \begin{align*}
        \frac{\norm{b-A\hatx_{i}}}{\norm{A}\norm{x}}&\leq \bigO\bigl(n(\kstar)^2\macheps\bigr)
        \max_{j\leq \kstar+1}\biggl(\frac{\norm{\hatx_j}}{\norm{x}}, 1\biggr)\qquad\text{and} \\
        \frac{\normA{\hatx_{i}-x}}{\norm{A}^{1/2}\norm{x}}
        &\leq \bigO\bigl(n(\kstar)^2\macheps\bigr) \kappa(A)^{1/2} \max_{j\leq \kstar+1}\biggl(\frac{\norm{\hatx_j}}{\norm{x}}, 1\biggr)
    \end{align*}
    with \(i\leq\kstar\), under the assumption
    \begin{equation*} 
    \begin{split}
        &\bigO\bigl(n\kstar\macheps\bigr)\kappa(A)
        + \bigO\bigl((\kstar)^2\macheps\bigr) \kappa(A)^{\frac{1}{2}}\leq \frac{1}{2}.
    \end{split}
    \end{equation*}
    In~\cite{Greenbaum1989}, the author implicitly provides an upper bound, i.e., \(\kstar\leq n^2\), with the restriction to achieve the forward error proportional to \(\sqrt{\macheps}\) instead of \(\macheps\).
\end{remark}

\subsection{Finite precision PCG using Cholesky factors as a left, right, or split preconditioner}
In this subsection, we discuss the backward and forward error results derived from Theorem~\ref{thm:backward-forward-err} for the PCG algorithm (i.e., Algorithm~\ref{alg:pcg}) using Cholesky factors where \(M=LL\trans\).
These factors can be employed as left (\(M_L=M\) and \(M_R=I\)), right (\(M_L=I\) and \(M_R=M\)), or split (\(M_L=L\) and \(M_R=L\trans\)) preconditioners.

We start by analyzing the rounding errors from applying the Cholesky factor to a vector in the form of~\eqref{eq:hatsk1} and~\eqref{eq:hatqk1}.

\begin{lemma} \label{lem:triangular-solver}
    Assume that \(L\in\mathbb{R}^{n\times n}\) is the Cholesky factor of \(M\) that satisfies \(M=LL\trans\) and let \(y_b\in \mathbb{R}^n\).
    If \(\bigO(n\macheps)\kappa(M)^{1/2}\leq 1/2\), then the result \(\hat y\) computed by solving \(Ly=y_b\) satisfies
    \begin{equation} \label{eq:lem:triangular-solver}
        \hat y = L^{-1}y_b + \Delta y\qquad\text{with}\qquad
        \norm{\Delta y}\leq \bigO(n\macheps)\kappa(M)^{\frac{1}{2}}\norm{L^{-1}y_b},
    \end{equation}
    and the result \(\hat y\) computed by solving \(L\trans y=y_b\) satisfies
    \begin{equation} \label{eq:lem:triangular-solver-trans}
        \hat y = L\itrans y_b + \Delta y\qquad\text{with}\qquad
        \norm{\Delta y}\leq \bigO(n\macheps)\kappa(M)^{\frac{1}{2}}\norm{L\itrans y_b}.
    \end{equation}
    Moreover, if \(\bigO(n\macheps)\kappa(M)\leq 1/2\), then the result \(\hat y\) computed by solving \(M y=y_b\) satisfies
    \begin{equation} \label{eq:lem:solveM}
        \hat y = M^{-1}y_b + \Delta y\qquad\text{with}\qquad
        \norm{\Delta y}\leq \bigO(n\macheps)\kappa(M) \norm{M^{-1} y_b}.
    \end{equation}
\end{lemma}

\begin{proof}
    From~\cite[Theorem 8.5]{H2002}, we have
    \begin{equation*}
        (L+\Delta L)\hat y = y_b\qquad\text{with}\qquad\abs{\Delta L}\leq \bigO(n\macheps)\abs{L},
    \end{equation*}
    which implies that
    \begin{equation} \label{eq:lem:proof:deltay-1}
    \begin{split}
        \norm{\hat{y} - L^{-1} y_b} \leq \bigO(n\macheps) \kappa(L)\norm{\hat{y}} 
        \leq \bigO(n\macheps) \kappa(L) (\norm{\hat{y} - L^{-1} y_b} + \norm{L^{-1} y_b}).
    \end{split}
    \end{equation}
    Together with the assumption \(\bigO(n\macheps)\kappa(M)^{1/2}\leq 1/2\), we have \(\bigO(n\macheps)\kappa(L)\leq 1/2\) and then
    \begin{equation} \label{eq:lem:proof:yhat-Linvyb}
    \begin{split}
        \norm{\Delta y} = \norm{\hat{y} - L^{-1} y_b}
        &\leq \frac{\bigO(n\macheps)\kappa(L)}{1-\bigO(n\macheps)\kappa(L)} \norm{L^{-1} y_b} \\
        &\leq \bigO(n\macheps)\kappa(L) \norm{L^{-1} y_b} \\
        &= \bigO(n\macheps)\kappa(M)^{1/2} \norm{L^{-1} y_b},
    \end{split}
    \end{equation}
    which yields~\eqref{eq:lem:triangular-solver}.
    The derivation of bound~\eqref{eq:lem:triangular-solver-trans} follows analogously from the proof of~\eqref{eq:lem:triangular-solver}.

    Next, we prove~\eqref{eq:lem:solveM}.
    Once again from~\cite[Theorem 8.5]{H2002}, \(\hat y\) satisfies
    \begin{equation*}
        (L+\Delta L)(L\trans + \Delta \tilde{L})\hat y = y_b\quad\text{with}\quad\abs{\Delta L}\leq \bigO(n\macheps)\abs{L}
        \quad\text{and}\quad\abs{\Delta \tilde{L}}\leq \bigO(n\macheps)\abs{L\trans},
    \end{equation*}
    which can be simplified as
    \begin{equation*}
        \bigl(M + \underbrace{\Delta LL\trans + L\Delta\tilde{L} + \Delta L\Delta\tilde{L}}_{=:\Delta M}\bigr)\hat y = y_b\quad\text{with}\quad\norm{\Delta M}\leq \bigO(n\macheps)\norm{L}^2 = \bigO(n\macheps)\norm{M}.
    \end{equation*}
    Similarly to~\eqref{eq:lem:proof:deltay-1} and~\eqref{eq:lem:proof:yhat-Linvyb}, we can conclude~\eqref{eq:lem:solveM}.
\end{proof}

Lemma~\ref{lem:triangular-solver} provides the values of \(\epss\) and \(\epsq\) for using Cholesky factors as left, right, and split preconditioners, respectively.

We now present sufficient conditions for the application of the preconditioners in different precisions such that the bounds from Theorem~\ref{thm:backward-forward-err} hold. In doing so, we distinguish between the three preconditioning schemes. Note that in the case of left and right preconditioning, both schemes culminate in the same Algorithm~\ref{alg:pcg} and possess the same theoretical result, as shown in the following Corollary.

\begin{corollary} \label{cor:mixed}
    Suppose that \(\hatxk\) and \(\hatrk\) are generated by Algorithm~\ref{alg:pcg} with \((M_L, M_R)=(M, I)\), \((M_L, M_R)=(I, M)\), or \((M_L, M_R) = (L, L^\top)\) where \(M=LL\trans\) is the Cholesky factorization of $M$.
    Moreover, assume that \(\machepss\) and \(\machepsq\), respectively, denote the unit roundoff values associated with the precisions involved in the application of the preconditioners \(M_L^{-1}\) and \(M_R^{-1}\) to a given vector in~\eqref{eq:hatsk1} and~\eqref{eq:hatqk1},
    and \(\barxk\) satisfies~\eqref{eq:tildexk1}.
    If \(\bigO(n\macheps)\kappa(A)+\bigO(n\macheps)\kappa(M)^{\frac{1}{2}}+\bigO(n\macheps_1)\kappa(M)^{m_1}\leq 1/2\), where
    \[
        \macheps_1 = 
        \begin{cases}
            \machepss, &\text{ for left PCG,}\\
            \machepsq, &\text{ for right PCG,}\\
            \machepss + \machepsq, &\text{ for split PCG,}\\
        \end{cases} \quad \text{and} \quad
        m_1 =
        \begin{cases}
             \frac{3}{2}, &\text{ for left or right PCG,}\\
             1, &\text{ for split PCG,}\\
        \end{cases}
    \]
    respectively, then there exists an iteration \(\kstar\geq 0\) such that
    \[
        f(\barx_{\kstar}) - f(\barx_{\kstar+1}) \leq c(\kstar)\cdot\bigO\bigl(n^2\macheps^2\bigr)\hatalpha_{\kstar}\norm{A}^2\norm{M^{-1}} \norm{x}^2,
    \]
    where \(c(\kstar)\) has been defined in Lemma~\ref{lem:Dfsmall}.
    Furthermore, if \(\kstar\) satisfies
    \begin{equation}  \label{eq:cor:assump-mixed}
    \begin{split}
        &\bigO\bigl(n\kstar\macheps\bigr)\kappa(A) + 
        \bigO\bigl(n\kstar\macheps\bigr)\kappa(M)^{\frac{1}{2}} + \bigO\bigl(n\kstar\macheps_1\bigr)\kappa(M)^{m_1} \\
        &+ \bigO\bigl((\kstar)^2\macheps) \kappa(M)^{\frac{1}{2}} \kappa\bigl(M_L^{-1}AM_R^{-1}\bigr)^{\frac{1}{2}}\leq \frac{1}{2},
    \end{split}
    \end{equation}
    then there exists \(i\leq \kstar\) such that
    \begin{align}
        \frac{\norm{\hatr_{i}}}{\norm{A}\norm{x}}&\leq \bigO\bigl(n(\kstar)^2\macheps\bigr) \kappa(M)^{\frac{1}{2}}
        \max_{j\leq \kstar+1}\biggl(\frac{\norm{\hatx_j}}{\norm{x}}, 1\biggr), \\
        \frac{\norm{b-A\hatx_{i}}}{\norm{A}\norm{x}}&\leq \bigO\bigl(n(\kstar)^2\macheps\bigr) \kappa(M)^{\frac{1}{2}}
        \max_{j\leq \kstar+1}\biggl(\frac{\norm{\hatx_j}}{\norm{x}}, 1\biggr),
    \end{align}
    and
    \begin{equation}
        \begin{split}
            \frac{\normA{\hatx_{i}-x}}{\norm{A}^{1/2}\norm{x}}
            &\leq \bigO\bigl(n(\kstar)^2\macheps\bigr)\kappa(M)^{\frac{1}{2}} \kappa(A)^{1/2} \max_{j\leq \kstar+1}\biggl(\frac{\norm{\hatx_j}}{\norm{x}}, 1\biggr).
        \end{split}
    \end{equation}
\end{corollary}

\begin{proof}
    We first consider the left preconditioned case, i.e, \(M_L=M\) and \(M_R=I\).
    By Lemma~\ref{lem:triangular-solver}, \(\epss=\bigO(n\machepss)\kappa(M)\) and \(\epsq=0\).
    Based on Theorem~\ref{thm:backward-forward-err}, we only need to estimate \(\epspre\), which requires the bound on \(\max_{y\neq 0} \bigl(\norm{M_R^{-1}y}\norm{M_L^{-1}y}/\normMinv{y}^2\bigr)\):
    \begin{equation} \label{eq:lem:proof:maxMRyMLy/yMnorm}
        \max_{y\neq 0} \biggl(\frac{\norm{M_R^{-1}y}\norm{M_L^{-1}y}}{\normMinv{y}^2}\biggr)
        \leq \max_{y\neq 0} \biggl(\frac{\norm{M^{-\frac{1}{2}}}\norm{M^{-\frac{1}{2}}y}\norm{M^{\frac{1}{2}}}\normMinv{y}}{\normMinv{y}^2}\biggr)
        \leq \kappa(M)^{\frac{1}{2}}.
    \end{equation}
    Together with the definition of \(\epspre\) in Lemma~\ref{lem:eachstep-alphabeta}, this implies that
    \begin{equation*}
        \epspre\leq \bigO\bigl(n\machepss\bigr)\kappa(M)^{\frac{3}{2}}
        +\bigO\bigl(n\macheps\bigr)\kappa(M)^{\frac{1}{2}},
    \end{equation*}
    which concludes the proof for the left preconditioned case.

    Similarly, for the right preconditioned case, i.e., \(M_L=I\) and \(M_R=M\), also by Lemma~\ref{lem:triangular-solver}, we have \(\epss=0\) and \(\epsq=\bigO(n\machepsq)\kappa(M)\).
    From the estimation of the quantity \(\max_{y\neq 0} \bigl(\norm{M_R^{-1}y}\norm{M_L^{-1}y}/\normMinv{y}^2\bigr)\) in~\eqref{eq:lem:proof:maxMRyMLy/yMnorm}, we have
    \begin{equation*}
        \epspre\leq \bigO\bigl(n\machepsq\bigr)\kappa(M)^{\frac{3}{2}}
        +\bigO\bigl(n\macheps\bigr)\kappa(M)^{\frac{1}{2}},
    \end{equation*}
    which concludes the proof for the right preconditioned case. 
    
    Finally, we prove the bounds for the split preconditioned case.
    We have \(M_L=L\) and \(M_R=L\trans\).
    By Lemma~\ref{lem:triangular-solver}, we have \(\epss=\bigO(n\machepss)\kappa(M)^{1/2}\) and \(\epsq=\bigO(n\machepsq)\kappa(M)^{1/2}\).
    Based on Theorem~\ref{thm:backward-forward-err}, we only need to estimate \(\epspre\), which requires the bound on \(\max_{y\neq 0} \bigl(\norm{M_R^{-1}y}\norm{M_L^{-1}y}/\normMinv{y}^2\bigr)\):
    \begin{equation*}
        \max_{y\neq 0} \biggl(\frac{\norm{M_R^{-1}y}\norm{M_L^{-1}y}}{\normMinv{y}^2}\biggr)
        = \max_{y\neq 0} \biggl(\frac{\norm{L^{-1}y}\norm{L\itrans y}}{\normMinv{y}^2}\biggr)
        = \max_{y\neq 0} \biggl(\frac{\sqrt{y\trans M^{-1}y}\sqrt{y\trans M\itrans y}}{\normMinv{y}^2}\biggr)
        = 1.
    \end{equation*}
    Together with the definition of \(\epspre\) in Lemma~\ref{lem:eachstep-alphabeta}, this implies that
    \begin{equation*}
        \epspre\leq \bigO\bigl(n(\machepss+\machepsq)\bigr)\kappa(M)
        +\bigO\bigl(n\macheps\bigr)\kappa(M)^{\frac{1}{2}},
    \end{equation*}
    which concludes the proof.
\end{proof}

\begin{remark}\label{remark:corollary}
    From Corollary~\ref{cor:mixed}, as long as the assumption in~\eqref{eq:cor:assump-mixed} is satisfied, which depends on the precisions used for applying the preconditioners, there is no effect on the backward and forward errors.
    When comparing the assumption on \(\kstar\) in the split preconditioned scenario with that for the left or right preconditioned cases, one notices that the second term in~\eqref{eq:cor:assump-mixed}, which directly affects the permissible low precision for applying the preconditioners, depends on \(\kappa(M)\) instead of \(\kappa(M)^{3/2}\).
    Consequently, in theory, PCG with split preconditioners enables the application of preconditioners at a lower precision level than its left- and right-preconditioned counterparts.
    Therefore, split-preconditioned PCG is preferred theoretically.
\end{remark}

\section{Numerical Experiments}
\label{sec:experiments}
In the following, we present a set of numerical experiments designed to illustrate the theoretical results in the previous section. These experiments incorporate preconditioning and application of the preconditioners in lower precisions. The implementation of our PCG framework can be found in the software package \texttt{MixedPrecisionPCG.jl}\footnote{\url{https://github.com/thbake/MixedPrecisionPCG.jl}} implemented in the Julia programming language. For simplicity, we slightly abuse the previously established notation and denote unit roundoffs associated with the application of left and right preconditioners by $\machepsleft$ and $\machepsright$, respectively. To be precise, we set $\machepsleft = \machepss$ and $\machepsright = \machepsq$. Moreover, we consider the arithmetics listed in Table~\ref{tab:arithmetics} for the storage and application of left and right preconditioners and set the working precision $\macheps = 1.11\times 10^{-16}$ (fp64 arithmetic) for every experiment. 
Apart from the application of the preconditioners, all computations---including the construction of matrices and preconditioners and the evaluation of errors---are all performed in the working (double) precision.
If we say that a preconditioning scheme was applied in a certain precision, we mean that the application of the preconditioner was carried out in said precision.

\begin{table}[!tb]
\centering
\caption{Floating point arithmetics used in the experiments and their associated unit roundoffs and smallest positive subnormal numbers.}
\label{tab:arithmetics}
\begin{tabular}{lcc}
    \hline
    Arithmetic & Unit roundoff & Smallest positive subnormal number\\
    \hline
    bfloat16      & $3.91 \times 10^{-3}$  & $9.18\times 10^{-41}$\\ 
    fp16 (half)   & $4.88 \times 10^{-4}$  & $5.96\times 10^{-8}$\\
    fp32 (single) & $5.96 \times 10^{-8}$  & $1.40\times 10^{-45}$\\
    fp64 (double) & $1.11 \times 10^{-16}$ & $4.94\times 10^{-324}$\\
    \hline
\end{tabular}
\end{table}

In order to generate linear algebraic systems and have control over the distribution of eigenvalues of the corresponding matrices, we follow the setting presented in \cite[Section 2.2]{CLS2024}. More precisely, given $0 < \lambda_1 < \lambda_n$ we consider the diagonal matrix $A = \diag(\lambda_1, \lambda_2, \dots, \lambda_n)$  with eigenvalues 
\begin{equation}\label{eq:exp:strakos-matrix}
    \lambda_i = \lambda_1 + \frac{i - 1}{n - 1} (\lambda_n - \lambda_1) \rho^{(n - i)}, \quad i = 2, \dots, n-1, \quad \rho \in [0,1].
\end{equation}
The right-hand side is given by $b = (1/\sqrt{n}) \cdot \bmat{1& \cdots& 1}\trans$ and we use a zero initial guess in all our experiments. The parameter $\rho$ in \eqref{eq:exp:strakos-matrix} controls the distribution of eigenvalues.
For each of the following experiments we consider $n = 85$ and set $\lambda_1 = 1, \lambda_n = 10^5,$ and $\rho = 0.6$, which yields $\kappa(A) = 10^5$, and an accumulation of eigenvalues on the left side of the spectrum.
We choose this distribution intentionally to have a convergence rate that is slow enough to benefit from the application of a preconditioner.
To construct the preconditioner, we use a similar setup to that of~\cite[Section 7.2.3.1]{V2022}.
That is, given an index $i \in \{2, \dots, n\}$, we employ 
\begin{equation}
    M = \diag(\lambda_1, \dots, \lambda_{i - 1}, \lambda_i, \dots, \lambda_i),
\end{equation} 
as the preconditioner
so that we truncate the last $n - i$ eigenvalues of \(A\).
In the case of split preconditioning, we take the Cholesky factors of the diagonal matrix \(M\) and apply it as both the left and right preconditioner, that is, \(M = M_L M_R = L L\trans\), computed using the \texttt{sqrt} command.
Then the application of preconditioners to vectors is computed by solving triangular systems.
In this way, the condition number of the preconditioned matrices $M^{-1} A$ and $A M^{-1}$ is given by $\kappa(M^{-1} A) = \kappa(A M^{-1}) = \lambda_n / \lambda_i < \lambda_n / \lambda_1 = \kappa(A)$ (in exact arithmetic). The spectrum of these matrices is given by
\[
	S = \left\{1, \frac{\lambda_{i + 1}}{\lambda_i}, \dots, \frac{\lambda_n}{\lambda_i} \right\},
\]
where the first $i$ eigenvalues are equal to one with multiplicity $i$.
The matrices $M^{-1} A$ and $L^{-1} A L^{-\top}$ are similar, and hence $S$ corresponds to both the spectrum of the left and split preconditioned matrices. 

By \eqref{eq:polynomial-eigvals-min} we know that the error in the $A$-norm is related to minimizing the squared values of a polynomial $p \in \mathcal{P}_k (0)$ at the eigenvalues of $A$. 
The mechanism of the CG algorithm to achieve this minimization is the placement of Ritz values close to the eigenvalues of $A$ in each iteration; see \cite[Section 5.6.1 and Section 5.9.1]{LS2012} for detailed expressions of the error in the $A$-norm involving Ritz values, and for an illuminating example that describes this mechanism, respectively.
In exact arithmetic with a sufficiently small value of $i$, we will have slow convergence due to the fact that the majority of terms in the sum in \eqref{eq:polynomial-eigvals-min} involve dividing by small eigenvalues of $A$.
Hence, the squared values of the CG polynomial $\varphi^\text{CG}_k$ -- which minimizes the right-hand side in \eqref{eq:polynomial-eigvals-min} -- at the eigenvalues of $A$ need to be small.
By the previous argument, the CG polynomial needs to place its roots (the Ritz values) close to the small eigenvalues of $A$.
This placement of Ritz values requires many iterations in comparison to CG applied to a matrix with an accumulation of eigenvalues to the right end of the spectrum with a few small outlying eigenvalues.
If we have a very large value of $\varphi^\text{CG}_k$ at the large outlying eigenvalue (say $\lambda_n$) of $A$, then $\varphi^\text{CG}_k$ must necessarily have a root close to said eigenvalue.
Moreover, if $\varphi^\text{CG}_k(\lambda_n)^2$ is large enough so that the division by $\lambda_n$ does not produce a small error, the CG algorithm will progressively place roots that are arbitrarily close to $\lambda_n$.
Since we have an accumulation of eigenvalues on the left end of the spectrum with large outliers far away from zero, small values of $i$ will not have a significant effect on the convergence of PCG compared to its unpreconditioned counterpart. As the value of $i$ increases, the eigenvalues indexed by \(i+1\), \(i+2\), \(\dots\), \(n\) of the original matrix are ``brought closer'' to the left end of the spectrum. As explained in Section~\ref{sec:algorithm}, left and right PCG are mathematically and numerically equivalent for a preconditioner given by a Cholesky factorization so that there is no need to run both variants. 

\begin{figure}[!tb]
    \begin{tabular}{ll}
        \includegraphics[width=\textwidth]{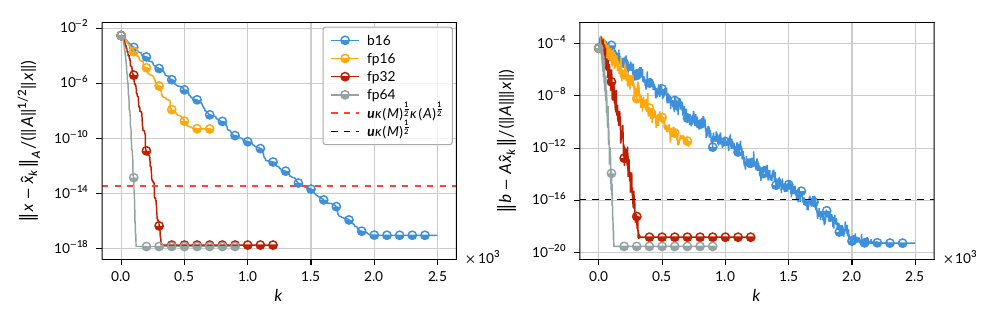}
    \end{tabular}
    \caption{Convergence curves for left PCG using the arithmetics listed in Table~\ref{tab:arithmetics} for the application of the preconditioners. Left: Left-hand side of \eqref{eq:thm-final:normtilder}. Right: Left-hand side of \eqref{eq:thm-final:normA-xk-x}.}
    \label{fig:leftPCG-bounds} 
\end{figure}

In our first experiment, we set $i = 55$ so that $\kappa(M^{-1} A) \approx 9.86 \times 10^4$ (in exact arithmetic), and test the bounds \eqref{eq:thm-final:normtilder} and \eqref{eq:thm-final:normA-xk-x} from Theorem~\ref{thm:backward-forward-err} for left PCG using arithmetics listed in Table~\ref{tab:arithmetics} by running the algorithm for 2500 iterations.
Figure~\ref{fig:leftPCG-bounds} displays the resulting convergence curves corresponding to the 2-norm of the true residual as well as the $A$-norm of the error, both divided by their corresponding normed quantities as described in equations \eqref{eq:thm-final:normtilder} and \eqref{eq:thm-final:normA-xk-x}, respectively. Following Remark~\ref{remark:kstar}, we \emph{exclude} the influence of the dimension size and of the number of iterations when plotting the respective bounds. 
As expected, the convergence rates decrease dramatically as the precisions associated to the arithmetics decrease. However, apart from the curve linked to the use of fp16, each curve reaches a level of accuracy close to $\macheps$. 
Since the bounds include the square roots of the condition number of both the preconditioner and the system matrix as multiplicative factors, the level reached by the curves is below the respective bounds.
Notice that this holds for the application of the preconditioners in bfloat16 arithmetic even though $n (\machepsleft + \machepsright) \kappa(M^{-1}) > 1/2$ so that the assumption \eqref{eq:cor:assump-mixed} for the left preconditioned case, does not hold for any $k \geq 2$. For $k \approx 700$, the first triangular solve (see Line 9 in Algorithm~\ref{alg:pcg}) performed in fp16 arithmetic yields a vector $\hat{s}_k = \hat{q}_k$ with entries close to the smallest positive subnormal number (see Table~\ref{tab:arithmetics}).
As a consequence, the inner product $\hat{r}_k\trans \hat{q_k} = 0$ for any $\hat{r}_k$ with sufficiently small entries such that the computation of $\hat{\beta}_{k + 1}$ satisfies 

\[ 
    \hat{\beta}_{k + 1} = 
    \frac{ (\hat{r}_{k+1} + \Delta \hat{r}_{k + 1})\trans \hat{q}_{k + 1} }{ (\hat{r}_k + \Delta \hat{r}_k)\trans \hat{q}_k}(1 + \delta \hat{\beta}_{k + 1}) =  \texttt{NaN}.
\]
Hence, the process fails completely before reaching the respective bounds and also convergence, which explains the behavior displayed by the fp16 curves in Figure~\ref{fig:leftPCG-bounds}.\footnote{We believe that this problem can be avoided by employing some scaling techniques.}

\begin{figure}[!tb]
    \begin{tabular}{l}
          \includegraphics[width=\textwidth]{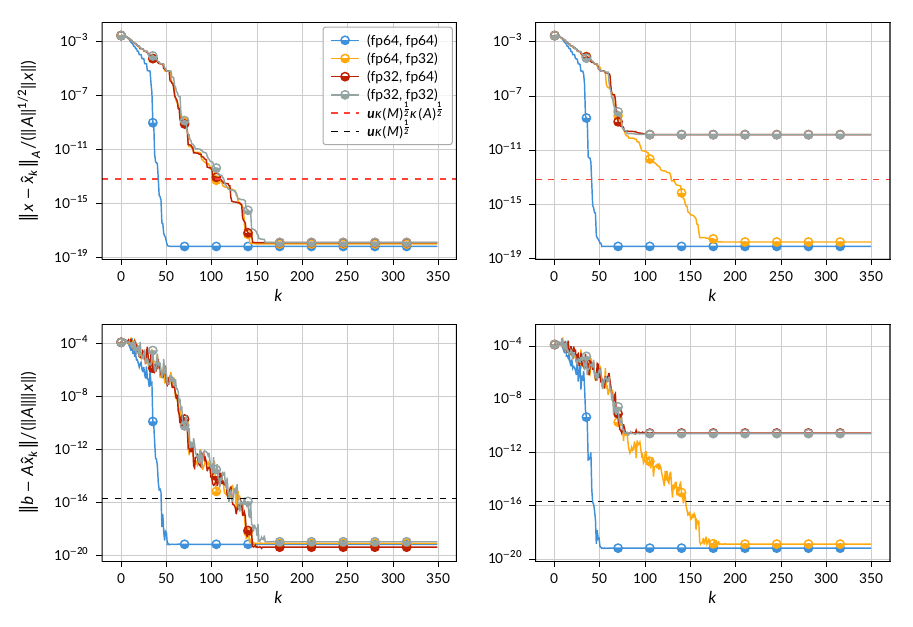}
    \end{tabular}
    \caption{Comparison of convergence curves for split PCG variants using combinations of double and single precisions for the application of the preconditioners. Left: Split PCG from Algorithm~\ref{alg:pcg}. Right: Saad's split PCG variant.}
    \label{fig:splitPCG-comparison} 
\end{figure}

For our next experiment, we set $i = 65$ to increase the rate of convergence, compare the accuracy reached by Algorithm~\ref{alg:pcg} and Saad's split PCG algorithm (see \cite[Algorithm 9.2]{Saad2003}) for all possible combinations of double and single precision with respect to the application of left and right preconditioners, and show the results in Figure~\ref{fig:splitPCG-comparison}. 
In the left column of Figure~\ref{fig:splitPCG-comparison} we observe that the PCG process of Algorithm~\ref{alg:pcg} reaches working precision accuracy, and it is therefore, below the respective bounds at every instance.
Furthermore, we observe that there is no significant difference in the convergence rate when at least one of the preconditioners is stored and applied in single precision.
However, this is not the case for \cite[Algorithm 9.2]{Saad2003} when the left preconditioner is applied in single precision, where both the \(2\)-norm of the true residual and the error in the $A$-norm eventually reach a point of stagnation as shown in the right column of Figure~\ref{fig:splitPCG-comparison}.
As described in Section~\ref{sec:algorithm}, this happens because the application of the left preconditioner is embedded into the computation of the recursively updated residual; that is, $r_{k + 1} = r_k - \alpha_k M_L^{-1} A p_k$ (see \cite[Line 5, Algorithm 9.2]{Saad2003}).
This in turn affects the norm of the residual gap, and it can be shown---following the analysis in the proof of Theorem~\ref{thm:Drk}---that the latter can be bounded by terms on the order of $\bigO(\macheps + \machepsleft)$.\footnote{In the interest of space, we exclude said derivation from this paper.
} 
Thus, whenever $\machepsleft$ is less than $\macheps$ the attainable accuracy of Saad's split PCG variant will be limited by the order of the precision with which the left preconditioner is applied.

\begin{figure}[!tb]
    \begin{tabular}{l l}
        \subcaptionbox{BE\label{fig:be-k2500}}[0.5\textwidth]
        {\includegraphics[trim=12 0 12 0, clip,width=0.5\textwidth]{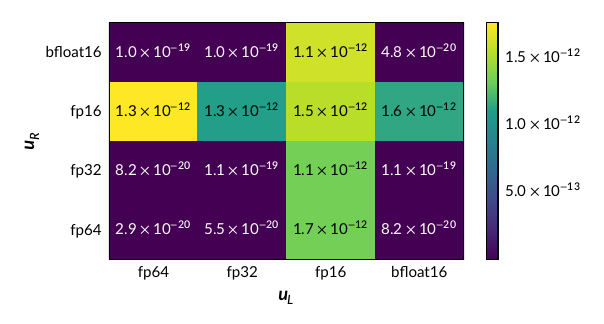}}
        &
        \subcaptionbox{FE\label{fig:fe-k2500}}[0.5\textwidth]
        {\includegraphics[trim=12 0 12 0, clip,width=0.5\textwidth]{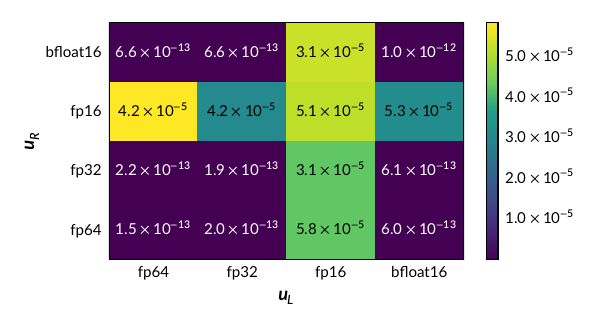}}\\
        \multicolumn{2}{c}{\subcaptionbox{IC\label{fig:ic-k2600}}[0.47\textwidth]{\includegraphics[width=0.50\textwidth]{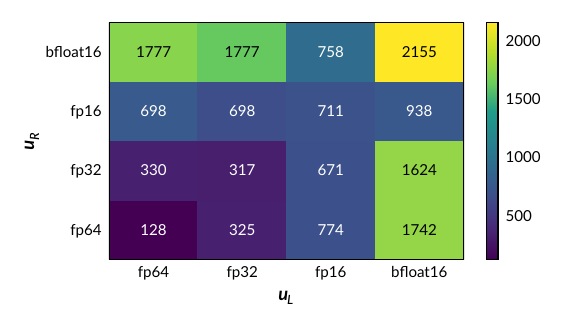}}}
    \end{tabular}
    \caption{Achieved backward and forward errors by split PCG for (first row) and corresponding number of iterations required (second row).}
    \label{fig:splitPCG-heatmaps} 
\end{figure}

\begin{figure}[!tb]
    \begin{tabular}{l}
        \includegraphics[width=\textwidth]{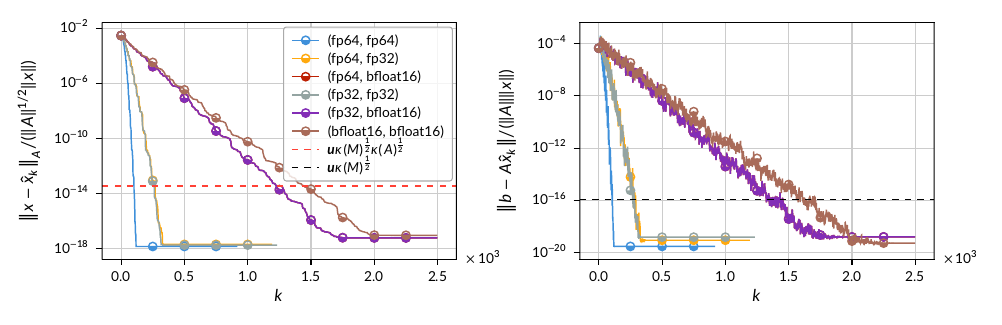}
    \end{tabular}
    \caption{Convergence curves for split PCG using arithmetics listed in Table~\ref{tab:arithmetics} excluding fp16.}
    \label{fig:splitPCG-bfloat16} 
\end{figure}

We now set $i = 55$ again and compute the same Cholesky factorization as in the first experiment.
Subsequently, we run split PCG for all 16 possible configurations using the arithmetics listed in Table~\ref{tab:arithmetics} until the error in the $A$-norm reaches its smallest value in each PCG process. Then, for the corresponding iteration steps $k$ we display the attained relative forward error (FE) in the $A$-norm, the relative backward error (BE), and the required iteration counts (IC) in the heatmaps in Figure~\ref{fig:splitPCG-heatmaps}. In the first row of Figure~\ref{fig:splitPCG-heatmaps} we can observe that for any configuration of the set $\{ \text{fp64, fp32, bfloat16} \}$ is close or surpasses the working precision. As expected from our previous observation of half precision, any configuration of the preconditioners involving fp16 arithmetic results in the failure to reach the desired accuracy due to underflow.
A closer look at the heatmaps reveals that, in general, there is no significant difference in the obtained accuracy whenever $\machepsleft \neq \machepsright$.
Figure~\ref{fig:splitPCG-bfloat16} shows the convergence curves for the split PCG setting, plotting the same metrics as in Figure~\ref{fig:leftPCG-bounds}\footnote{We exclude the curves involving fp16 arithmetic for visibility.
Essentially, we have the same behavior as in Figure~\ref{fig:leftPCG-bounds} and this case was thoroughly discussed.}.
Comparing these two figures, we see that there is no significant difference in the convergence behavior of left and split PCG.
Once again, for split PCG using bfloat16 arithmetic for the application of the left or right preconditioner, assumption \eqref{eq:cor:assump-mixed} is not met for any $k \geq 1$, but the method still converges.
These last two points suggest that the discussion of Remark~\ref{remark:corollary} seems more theoretical than of practical relevance.

%
%
As described in \cite[Section 2]{CLS2024}, finite precision CG applied to the linear system with $A$ that has the eigenvalue distribution described above will showcase a significant delay of convergence. The reason for this delay is the placement of \emph{multiple} Ritz values to decrease the otherwise huge gradients of the CG polynomial, since the proximity with which roots are placed in finite precision computations is limited by the machine precision.
As observed in our experiment, this delay of convergence is exacerbated with increasing unit roundoff size associated with the application of the preconditioners.
Thus, if one is willing to use low precision to apply the preconditioners---e.g., if there is some potential performance gain despite the delay of convergence---it seems that the most reasonable choice is to take $\machepsleft = \machepsright$.

To conclude the set of experiments, we adapt a recent synthetic example from \cite[Section 4]{EGN2025} that displays deterioration in the attainable accuracy with a preconditioner, even with a zero initial guess.
Here, the problem size is \(n = 100\), and we have a dense matrix \(A = U \Sigma U\trans \).
The diagonal matrix \( \Sigma \) has logarithmically spaced eigenvalues between \(1 \text{ and } 10^{12}\) resulting in \(\kappa(A) = 10^{12}\), and \(U \in \mathbb R^{n \times n} \) is a Haar random orthogonal matrix. 
We employ an SPD preconditioner \(P\) that satisfies \(\kappa (P) \approx 10^{12} \text{ and } \kappa(P^{-1/2} A P^{-1/2}) \approx 7.8 \).
Its construction involves the matrices \(U, \Sigma \) and a Wishart matrix; for more details on the construction of this example, we refer to \cite[Appendix A]{EGN2025}.

\begin{figure}[!tb]
    \begin{tabular}{l}
          \includegraphics[width=\textwidth]{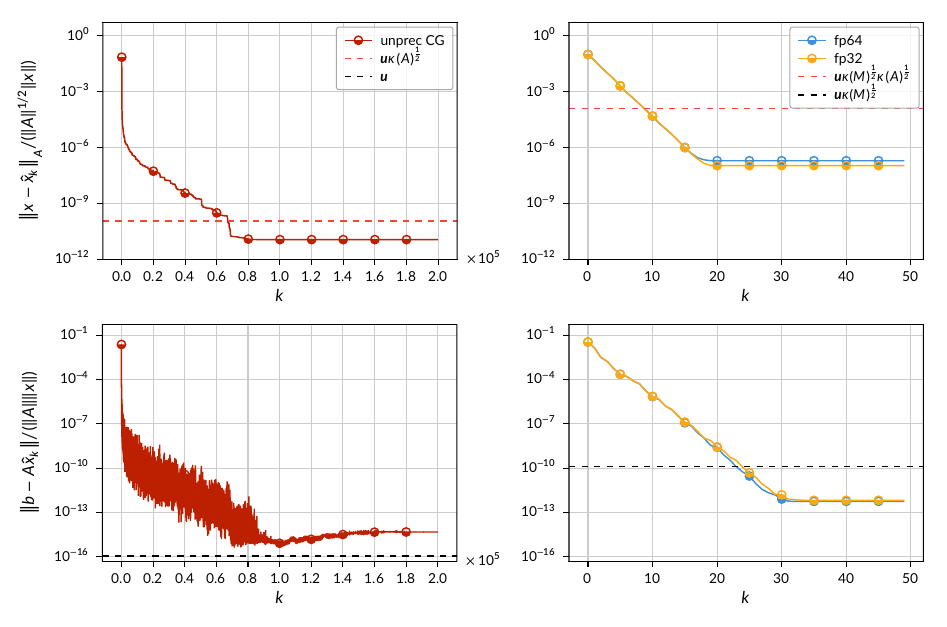}
    \end{tabular}
    \caption{Convergence curves for unpreconditioned CG (left column) and left PCG (right column) applied to the first problem in \cite[Section 3]{EGN2025} with double and single precision preconditioner application.}
    \label{fig:bad-example} 
\end{figure}

We run unpreconditioned CG, and left PCG applying the Cholesky factors of \(P\) in double and single precision, and plot the same accuracy metrics as in our previous examples in Figure~\ref{fig:bad-example}.
For unpreconditioned CG our bounds are particularly descriptive, since both the scaled \(A\)-norm of the error and the scaled true residual norm reach (or are very close) to \(\macheps \kappa(A)^{1/2}\) and \(\macheps\), respectively.
In contrast, PCG stagnates orders of magnitude before reaching these quantities, which shows the influence of the preconditioner on the attainable accuracy.\footnote{The authors in \cite{EGN2025} show experimentally that this loss of accuracy can be overcome with (a few steps of) iterative refinement.}
Despite a more ill-conditioned matrix and preconditioner than in the previous examples, the condition number alone does not seem to fully explain the loss of attainable accuracy since the same example without using preconditioners does converge to unit roundoff accuracy.
Last, it is worth emphasizing that this stagnation is not caused by use of lower precision either, since both instances of double and single precision stagnate at the same level. 

\section{Conclusions}
\label{sec:conclusions}
Current finite precision analyses of the CG algorithm mainly focus on analyzing the residual gap based on the assumption that the norm of the recursively updated residual \(\hatrk\) goes to a sufficiently small level after enough iterations.
For a long time, a rigorous proof has been lacking to confirm that the norm of \(\hatrk\) can achieve a level close to the machine precision after sufficiently many iterations.
In this work, we present a rigorous analysis to show that \(\norm{\hatrk}/(\norm{A}\norm{x})\) can reach the level \(\bigO(n(\kstar)^2\macheps)\) at iteration \(k = \kstar\), provided that a certain moderate assumption on the specific iteration \(\kstar\) is met. 
Combining this result with a bound on the residual gap, we further demonstrate that PCG can achieve a backward and forward stable solution before  iteration \(\kstar\).
Our analysis demonstrates the worst-case scenario, highlighting that the presence of \(n\) and \(\kstar\) in the assumption renders it challenging to satisfy, particularly for large linear systems.
However, in practical applications, the influence of \(n\) and \(\kstar\) on rounding errors is generally much less significant before convergence.
Once the required accuracy is achieved, the presence of the iteration number in the bounds suggests that an appropriate stopping criterion is crucial to cease PCG iterations.
This prevents the accumulation of rounding errors as \(k\) continues to increase, which might lead to an increase in the residual after it decreases to \(\bigO(\macheps)\).
In addition, our analysis also ensures that one can employ lower precision in applying preconditioners while still achieving working precision accuracy for the computed approximation, as long as reasonable conditions are satisfied, which is verified by numerical experiments.
Finally, we have introduced a PCG framework that encompasses all left, right, and split preconditioning schemes. 
In contrast to the classical split PCG algorithm introduced in \cite{Saad2003}, the attainable accuracy of our split preconditioned variant is not affected when low precision is employed for the application of the left preconditioner.
\section*{Acknowledgments}

All authors are supported by the European Union (ERC, inEXASCALE, 101075632). Views and opinions expressed are those of the authors only and do not necessarily reflect those of the European Union or the European Research Council. Neither the European Union nor the granting authority can be held responsible for them. E.~Carson and Y.~Ma additionally acknowledge support from the Charles University Research Centre program No. UNCE/24/SCI/005.

\bibliographystyle{abbrvurl}
\bibliography{mybib}


\end{document}